\documentclass[10.5pt]{article}
\usepackage{amsfonts}

\def\sqr#1#2{{\vcenter{\vbox{\hrule height.#2pt
              \hbox{\vrule width.#2pt height#1pt \kern#1pt \vrule width.#2pt}
              \hrule height.#2pt}}}}
\def\signed #1{{\unskip\nobreak\hfil\penalty50
              \hskip2em\hbox{}\nobreak\hfil#1
              \parfillskip=0pt \finalhyphendemerits=0 \par}}
\def\endpf{\signed {$\sqr69$}}
\def\3n{\negthinspace \negthinspace \negthinspace }
\def\2n{\negthinspace \negthinspace }
\def\1n{\negthinspace }

\def\dbE{\mathbb{E}}
\def\dbF{\mathbb{F}}

\def\dbP{\mathbb{P}}

\def\dbR{\mathbb{R}}

\def\={\buildrel \triangle \over =}

\def\ds{\displaystyle}

\def\ns{\noalign{\ss}}

\def\b{\beta}
\def\g{\gamma}
\def\d{\delta}

\def\z{\zeta}

\def\l{\lambda}
\def\m{\mu}

\def\si{\sigma}
\def\t{\tau}
\def\f{\varphi}

\def\o{\omega}

\def\D{\Delta}

\def\F{\Phi}
\def\O{\Omega}

\def\cA{{\cal A}}

\def\cF{{\cal F}}

\def\cH{{\cal H}}

\def\cM{{\cal M}}

\def\cl{{\cal l}}
\def\ss{\smallskip}
\def\ms{\medskip}

\def\q{\quad}
\def\qq{\qquad}
\def\hb{\hbox}

\def\lan{\mathop{\langle}}
\def\ran{\mathop{\rangle}}

\def\esssup{\mathop{\rm esssup}}

\def\wt{\widetilde}
\def\cd{\cdot}
\def\cds{\cdots}

\def\ae{\hbox{\rm a.e.{ }}}
\def\as{\hbox{\rm a.s.{ }}}

\def\cl{\overline}

\def\({\Big (}
\def\){\Big )}
\def\[{\Big[}
\def\]{\Big]}
\def\bde{\begin{definition}}
\def\ede{\end{definition}}
\def\be{\begin{equation}}
\def\bel{\begin{equation}\label}
\def\ee{\end{equation}}
\def\bt{\begin{theorem}}
\def\et{\end{theorem}}
\def\bc{\begin{corollary}}
\def\ec{\end{corollary}}
\def\bl{\begin{lemma}}
\def\el{\end{lemma}}
\def\bp{\begin{proposition}}
\def\ep{\end{proposition}}
\def\bas{\begin{assumption}}
\def\eas{\end{assumption}}
\def\br{\begin{remark}}
\def\er{\end{remark}}
\def\ba{\begin{array}}
\def\ea{\end{array}}
\def\ed{\end{document}}

\def\square#1{\vbox{\hrule\hbox{\vrule height#1%
     \kern#1\vrule}\hrule}}
\def\rectangle#1#2{\vbox{\hrule\hbox{\vrule height#1%
     \kern#2\vrule}\hrule}}

\font\tenbb=msbm10 \font\sevenbb=msbm7 \font\fivebb=msbm5

\newfam\bbfam
\scriptscriptfont\bbfam=\fivebb \textfont\bbfam=\tenbb
\scriptfont\bbfam=\sevenbb

\newtheorem{lemma}{Lemma}[section]
\newtheorem{remark}{Remark}[section]

\newtheorem{theorem}{Theorem}[section]
\newtheorem{corollary}{Corollary}[section]

\newtheorem{definition}{Definition}[section]
\newtheorem{proposition}{Proposition}[section]
\newtheorem{assumption}{Assumption}[section]

\makeatletter
   
   \@addtoreset{equation}{section}
\makeatother

\begin{document}
\title{\bf Comparison Theorems for Backward Stochastic Volterra Integral Equations\footnote{This work is supported
in part by National Natural Science Foundation of China (Grants
10771122 and 11071145), Independent Innovation Foundation of
Shandong University (Grant 2010JQ010), Graduate Independent
Innovation Foundation of Shandong University (GIIFSDU), and the NSF
grant DMS-1007514.}}

\author{Tianxiao Wang\footnote{School of Mathematics, Shandong
University, Jinan 250100, China},~~and~~Jiongmin
Yong\footnote{Department of Mathematics, University of Central
Florida, Orlando, FL 32816, USA.}}

\maketitle

\begin{abstract}
For backward stochastic Volterra integral equations (BSVIEs) in
multi-dimensional Euclidean spaces, comparison theorems are
established in a systematic way for the adapted solutions and
adapted M-solutions. For completeness, comparison theorems for
(forward) stochastic differential equations, backward stochastic
differential equations, and (forward) stochastic Volterra integral
equations (FSVIEs) are also presented. Duality principles are used
in some relevant proofs. Also, it is found that certain kind of
monotonicity conditions play crucial roles to guarantee the
comparison theorems for FSVIEs and BSVIEs to be true. Various
counterexamples show that the assumed conditions are almost
necessary in some sense.

\end{abstract}

\ms


\ms


\ms

\section{Introduction.}\label{1}

Throughout this paper, we let $(\O,\cF,\dbF,\dbP)$ be a complete
filtered probability space on which a one-dimensional standard
Brownian motion $W(\cd)$ is defined with $\dbF=\{\cF_t\}_{t\ge0}$
being its natural filtration augmented by all the $\dbP$-null sets.
We consider the following equation in $\dbR^n$, the usual
$n$-dimensional real Euclidean space:
\bel{BSVIE1}Y(t)=\psi(t)+\int_t^Tg(t,s,Y(s),Z(t,s),Z(s,t))ds
-\int_t^TZ(t,s)dW(s),\q t\in[0,T],\ee
which is called a {\it backward stochastic Volterra integral
equation} (BSVIE, for short). Such kind of equations have been
investigated in the recent years (see \cite{Lin 2002,Yong 2006,Yong
2007,Yong 2008,Wang-Shi 2010, Anh-Grecksch-Yong 2011} and references
cited therein). BSVIEs are natural extensions of by now
well-understood backward stochastic differential equations (BSDEs,
for short) whose integral form is as follows:
\bel{1.2}Y(t)=\xi+\int_t^Tg(s,Y(s),Z(s))ds-\int_t^TZ(s)dW(s),\q
t\in[0,T].\ee
See \cite{Pardoux-Peng 1990, El Karoui-Peng-Quenez 1997, Ma-Yong
1999, Yong-Zhou 1999} for some standard results on BSDEs. An
interesting result of BSDEs is the comparison theorem for the
adapted solutions. A little precisely, say, for $n=1$, if
$(Y^i(\cd),Z^i(\cd))$ is the adapted solution to the BSDE
(\ref{1.2}) with $(\xi,g(\cd))$ replaced by $(\xi^i,g^i(\cd))$
($i=0,1$) such that
\bel{xi0<xi1}\left\{\ba{ll}
\ns\ds\xi^0\le\xi^1,\qq\as,\\
\ns\ds g^0(t,y,z)\le
g^1(t,y,z),\qq\forall(t,y,z)\in[0,T]\times\dbR\times\dbR,~\as,\ea\right.\ee
then
\bel{Y0<Y1}Y^0(t)\le Y^1(t),\qq t\in[0,T],~\as\ee
The comparison theorem also holds for multi-dimensional BSDEs. We
refer the readers to \cite{Hu-Peng 2006} for details. Because of the
comparison theorem, one can use the adapted solutions to BSDEs as
dynamic risk measures or stochastic differential utility for
(static) random variables which could be the payoff of a European
contingent claim at the maturity.

\ms

Now, for BSVIEs, from mathematical point of view, it is natural to
ask if a comparison theorem similar to that for BSDEs hold for
solutions to BSVIEs. More precisely, if $(Y^i(\cd),Z^i(\cd\,,\cd))$
is the solution to BSVIE (\ref{BSVIE1}), in a proper sense, with
$(\psi(\cd),g(\cd))$ replaced by $(\psi^i(\cd),g^i(\cd))$, $i=0,1$,
and
\bel{1.5}\left\{\ba{ll}
\ns\ds\psi^0(t)\le\psi^1(t),\qq t\in[0,T],~\as,\\
\ns\ds g^0(t,s,y,z,\z)\le g^1(t,s,y,z,\z),\qq0\le t\le s\le
T,~y,z,\z\in\dbR,~\as\ea\right.\ee
Can we have the comparison relation (\ref{Y0<Y1})?

\ms

On the other hand, similar to BSDEs, if proper comparison theorems
hold for BSVIEs, then there will be some interesting applications of
BSVIEs in risk management and optimal investment/comsuption
problems. Let us elaborate in a little details.

\ms

It is common that in order to expect some returns from various
existing risky assets, one should hold them for possibly different
length of time period. The value of the positions for these assets
at some future time form a (not necessarily adapted) stochastic
process, for which people would like to measure the dynamic risks. A
simple illustrative example can be found in \cite{Yong 2007}. We
emphasize that the processes (not just random variables) for which
one wants to measure the dynamic risk are not necessarily adapted.
Dynamic risk measures for discrete-time processes have been
considered in the literature, see, for examples, \cite{FS2005,
PFM2006EP, AFP2012} and so on. On the other hand, static risk
measures for continuous-time processes were studied in
\cite{PFM2004, PFM2006FS}. We believe that BSVIEs should be a useful
tool in studying dynamic risk measures for (not necessarily adapted)
stochastic processes. Therefore, to establish comparison theorems
for BSVIEs becomes quite necessary.

\ms

The second relevant motivation comes from the study of general yet
realistic stochastic utility problem. The stochastic differential
utility was introduced and studied in \cite{DE1992, LQ2003}, where
the intertemporal consistency and Bellman's principle of optimality
is applicable. However, real problems are usually of
time-inconsistent nature. In fact, many experimental study on time
preference shows that the standard assumption of time consistency is
unrealistic. Moreover, substantial evidence also suggest that agents
are impatient about choices in the short term but are patient among
the long-term alternatives. Recently, some people are interested in
the following type of stochastic utility function
$$Y(t)=\dbE\[\int_t^T\ell(t,s)u(c(s))ds\big|\cF_t\],\qq t\in[0,T],$$
with $\ell(t,s)$ being the discount factor, see \cite{MN2010,
EMP2012, Yong2012-2}. We expect that comparison theorems of BSVIEs
will play an important role in formulating general stochastic
utility functions and investigating their properties such as
comparative risk aversion, risk aversion, etc., which will
substantially extend the results in \cite{DE1992}.

\ms

We will present applications of comparison theorems of BSVIEs in
finance and other related area in our future publications.

\ms

Now, returning to comparison theorems for BSVIEs, we point out that
unlike BSDEs, (\ref{1.5}) is not enough to ensure comparison
relation (\ref{Y0<Y1}), in general. Various counterexamples will be
presented. Due to the complicated situation for BSVIEs, the theory
of comparison for solutions to BSVIEs is much more richer than that
for BSDEs. The main purpose of this paper is to establish various
comparison theorems for solutions to BSVIEs in multi-dimensional
Euclidean spaces. To this end, we first will consider BSVIE
(\ref{BSVIE1}) with the generator $g(\cd)$ independent of $Z(s,t)$.
For such a case, in order the comparison theorem holds, one needs
some kind of monotonicity for the generator $g(\cd)$ and/or the free
term $\psi(\cd)$. Some examples will show that the conditions we
impose are almost necessary. The second case to be considered is
that the generator $g(\cd)$ depends on $Z(s,t)$ and independent of
$Z(t,s)$. For such a case, we are comparing adapted M-solution for
(\ref{BSVIE1}) introduced in \cite{Yong 2008}. It turns out that
under proper monotonicity conditions, we are able to obtain a
comparison theorem for adapted M-solutions, which is weaker than
that for the first case. More precisely, instead of (\ref{Y0<Y1}),
we can only have
$$\dbE\[\int_t^TY^0(s)ds\big|\cF_t\]\le\dbE\[\int_t^TY^1(s)ds\big|\cF_t\],
\qq t\in[0,T],~\as$$
This result corrects a relevant result in \cite{Yong 2006, Yong
2007}. Finally, inspired by \cite{BLP2009} and \cite{BMZ2012}, we
introduce a new notion, called {\it conditional $h$-solutions} for
BSVIEs (\ref{BSVIE1}), and briefly discuss the corresponding
comparison theorem by following similar ideas for the first two
cases.

\ms

Note that the proofs of above results are closely connected with the
comparison theorems of (forward) stochastic differential equations
(FSDEs, for short), (forward) stochastic Volterra integral equations
(FSVIEs, for short), and BSDEs (allowing the dimension $n>1$). For
completeness, we will present/recall some relevant results here.
Interestingly, even for FSDEs and BSDEs, our proofs are different
from those in \cite{Geib-Manthey 1994,PZ2006,Hu-Peng 2006},
respectively, and more straightforward.

\ms

The rest of the paper is organized as follows: In Section 2, we
present some comparison theorems of FSDEs, BSDEs and FSVIEs. In
Section 3, we establish several comparison theorems for BSVIEs from
three different perspectives. Various persuasive examples will be
presented to illustrate the obtained results. Finally, some
concluding remarks are collected in Section 4.

\section{Comparison theorem for FSDEs, FSVIEs, and BSDEs}

\rm

In this section, we are going to present comparison theorems for
FSDEs, FSVIEs, and BSDEs, allowing the dimension $n>1$. Some of them
are known. But our proofs are a little different.

\ms

Let us first make some preliminaries. Denote
$$\dbR^n_+=\Big\{(x_1,\cds,x_n)\in\dbR^n\bigm|x_i\ge0,~1\le i\le n\Big\}.$$
When $x\in\dbR^n_+$, we also denote it by $x\ge0$, and say that $x$
is {\it nonnegative}. By $x\le0$ and $x\ge y$ (for $x,y\in\dbR^n$),
we mean $-x\ge0$ and $x-y\ge0$, respectively. In what follows, we
let $e_i\in\dbR^n_+$ be the vector that the $i$-th entry is 1 and
all other entries are zero. Let
$$\left\{\ba{ll}
\ns\ds\dbR^{n\times m}_+=\Big\{A=(a_{ij})\in\dbR^{n\times
m}\bigm|a_{ij}\ge0,~1\le i\le n,~1\le j\le m\Big\},\\
\ns\ds\dbR^{n\times n}_{*+}=\Big\{A=(a_{ij})\in\dbR^{n\times
n}\bigm|a_{ij}\ge0,~i\ne j\Big\}\equiv \Big\{A\in\dbR^{n\times
n}\bigm|\lan Ae_i,e_j\ran\ge0,~i\ne j\Big\},\\
\ns\ds\dbR^{n\times n}_d=\Big\{A=(a_{ij})\in\dbR^{n\times
n}\bigm|a_{ij}=0,~i\ne j\Big\}\equiv \Big\{A\in\dbR^{n\times
n}\bigm|\lan Ae_i,e_j\ran=0,~i\ne j\Big\}.\ea\right.$$
Note that $\dbR^{n\times m}_+$ is the set of all $(n\times m)$
matrices with all the entries being nonnegative, $\dbR^{n\times
n}_{*+}$ is the set of all $(n\times n)$ matrices with all the
off-diagonal entries being nonnegative (no conditions are imposed on
the diagonal entries), and $\dbR^{n\times n}_d$ is the set of all
$(n\times n)$ diagonal matrices, with the diagonal entries allowing
to be any real numbers. Clearly, $\dbR^{n\times m}_+$ and
$\dbR^{n\times n}_{*+}$ are closed convex cones of $\dbR^{n\times
m}$ and $\dbR^{n\times n}$, respectively; $\dbR^{n\times n}_+$ is a
proper subset of $\dbR^{n\times n}_{*+}$; and $\dbR^{n\times n}_d$
is a proper subspace of $\dbR^{n\times n}$, contained in
$\dbR^{n\times n}_{*+}$. Also,
$$\dbR^{n\times n}_{*+}=\dbR^{n\times n}_++\dbR^{n\times n}_d\equiv\Big\{A+B\bigm|
A\in\dbR^{n\times n}_+,\,B\in\dbR^{n\times n}_d\Big\}.$$
Further, for $n=m=1$, one has
\bel{2.1}\dbR^{1\times1}_{*+}=\dbR^{1\times1}_d=\dbR,\qq\dbR^{1\times1}_+=\dbR_+\equiv[0,\infty).\ee
We have the following simple result whose proof is obvious.

\ms

\bf Proposition 2.1. \sl Let $A\in\dbR^{n\times m}$. Then
$A\in\dbR^{n\times m}_+$ if and only if
\bel{}Ax\ge0,\qq\forall x\ge0.\ee

\rm

\ms

Next, we introduce some spaces. Let $H=\dbR^n,\dbR^{n\times m}$,
etc. with $|\cd|$ beng its norm. For $1\le p,q<\infty$ and $0\le
s<t\le T$, define
$$\ba{ll}
\ns\ds L^p_{\cF_t}(\O;H)=\Big\{\xi:\O\to H\bigm|\xi\hb{ is
$\cF_t$-measurable, }\dbE|\xi|^p<\infty\Big\},\\
\ns\ds L^p_{\cF_T}(\O;L^q(s,t;H))=\Big\{X:[s,t]\times\O\to
H\bigm|X(\cd)\hb{ is
$\cF_T$-measurable},~\dbE\(\int_s^t|X(r)|^qdr\)^{p\over
q}<\infty\Big\},\\
\ns\ds L^p_{\cF_T}(\O;C([s,t];H))=\Big\{X:[s,t]\times\O\to
H\bigm|X(\cd)\hb{ is $\cF_T$-measurable, has continuous paths,
}\\
\ns\ds\qq\qq\qq\qq\qq\qq\qq\qq\qq\qq\qq\dbE\(\sup_{r\in[s,t]}|X(r)|^p\)<\infty\Big\},\\
\ns\ds L^q_{\cF_T}(s,t;L^p(\O;H))=\Big\{X:[s,t]\times\O\to
H\bigm|X(\cd)\hb{ is
$\cF_T$-measurable},~\int_s^t\(\dbE|X(r)|^p\)^{q\over
p}dr<\infty\Big\},\\
\ns\ds C_{\cF_T}([s,t];L^p(\O;H))\1n=\1n\Big\{X\1n:\1n[s,t]\to
L^p_{\cF_T}(\O;H)\bigm|X(\cd)\hb{ is continuous,
}\sup_{r\in[s,t]}\dbE|X(r)|^p<\infty\Big\}.\ea$$
The spaces with the above $p$ and/or $q$ replaced by $\infty$ can be
defined in an obvious way. Also, we define
$$L^p_\dbF(\O;L^q(s,t;H))=\Big\{X(\cd)\in L^p_{\cF_T}(\O;L^q(s,t;H))
\bigm|X(\cd)\hb{ is $\dbF$-adapted}\Big\}.$$
The spaces $L^p_\dbF(\O;C([s,t];H))$, $L^q_\dbF(s,t;L^p(\O;H))$, and
$C_\dbF([s,t];L^p(\O;H))$ (with $1\le p,q\le\infty$) can be defined
in the same way. For simplicity, we denote
$$L^p_\dbF(s,t;H)=L^p_\dbF(\O;L^p(s,t;H))=L^p_\dbF(s,t;L^p(\O;H)),\qq1\le p\le\infty.$$
Further, we denote
$$\D=\Big\{(t,s)\in[0,T]^2\Bigm|t\le s\Big\},\q\D^*=\Big\{(t,s)\in[0,T]^2\Bigm|t\ge s\Big\}\equiv
\cl{\D^c},$$
and let
$$\ba{ll}
\ns\ds L^p_\dbF(\D;H)=\Big\{Z:\D\times\O\to H\bigm|s\mapsto
Z(t,s)\hb{ is $\dbF$-progressively measurable on $[t,T]$, $\forall t\in[0,T]$,}\\
\ns\ds\qq\qq\qq\qq\qq\qq\qq\qq\int_0^T\dbE\(\int_t^T|Z(t,s)|^2ds\)^{p\over2}dt<\infty\Big\},\ea$$
$$\ba{ll}
\ns\ds L^p(0,T;L^2_\dbF(0,T;H))=\Big\{Z:[0,T]^2\times\O\to
H\bigm|s\mapsto
Z(t,s)\hb{ is $\dbF$-progressively measurable}\\
\ns\ds\qq\qq\qq\qq\qq\qq\qq\qq\hb{ on $[t,T]$, $\forall t\in[0,T]$,}
\int_0^T\dbE\(\int_0^T|Z(t,s)|^2ds\)^{p\over2}dt<\infty\Big\}.\ea$$
The spaces $L^\infty_\dbF(\D;H)$ and $L^\infty(0,T;L^2_\dbF(0,T;H))$
can be defined similarly. Then we denote
$$\left\{\ba{ll}
\ns\ds\cH^p_\D[0,T]=L^p_\dbF(0,T)\times L^p_\dbF(\D;\dbR^n),\\
\ns\ds\cH^p[0,T]=L^p_\dbF(0,T)\times
L^p(0,T;L^2_\dbF(0,T;\dbR^n)),\\
\ns\ds\cM^p[0,T]=\Big\{(y(\cd),z(\cd\,,\cd))\in\cH^p[0,T]\bigm|y(t)=\dbE
y(t)+\int_0^tz(t,s)dW(s),\;t\in[0,T]\Big\}.\ea\right.$$

\subsection{Comparison of solutions to FSDEs.}

For any $(s,x)\in[0,T)\times\dbR^n$, let us first consider the
following linear FSDE:
\bel{LSDE1}\left\{\ba{ll}
\ns\ds dX(t)=\(A_0(t)X(t)+b(t)\)dt+A_1(t)X(t)dW(t),\q
t\in[s,T],\\
\ns\ds X(s)=x,\ea\right.\ee
with $A_0(\cd)$ and $A_1(\cd)$ satisfying the following
assumption.

\ms

{\bf(FD1)} The maps $A_0(\cd),A_1(\cd)\in
L^\infty_\dbF(\O;C([0,T];\dbR^{n\times n}))$.

\ms

We point out here that if the diffusion in (\ref{LSDE1}) is replaced
by $A_1(t)X(t)+\si(t)$ for some $\si(\cd)\ne0$, then comparison
theorem might fail in general. Therefore, we restrict ourselves to
the above form. It is standard that under (FD1), for any
$(s,x)\in[0,T)\times\dbR^n$, $b(\cd)\in
L^2_\dbF(\O;L^1(s,T;\dbR^n))$, FSDE (\ref{LSDE1}) admits a unique
solution $X(\cd)\equiv X(\cd\,;s,x,b(\cd))\in L^2_\dbF(\O;C([s,T];$
$\dbR^n))$, and the following estimate holds:
\bel{}\dbE\[\sup_{t\in[s,T]}|X(t)|^2\]\le
K\Big\{|x|^2+\dbE\(\int_s^T|b(t)|dt\)^2\Big\}.\ee
Hereafter, $K>0$ represents a generic constant which can be
different from line to line. Let $\F(\cd\,,\cd)$ be the {\it
stochastic fundamental matrix} of $\{A_0(\cd),A_1(\cd)\}$, i.e.,
\bel{}\left\{\ba{ll}
\ns\ds d\F(t,s)=A_0(t)\F(t,s)dt+A_1(t)\F(t,s)dW(t),\qq t\in[s,T],\\
\ns\ds\F(s,s)=I.\ea\right.\ee
Then one has the following {\it variation of constant formula}:
\bel{variation of constant}\ba{ll}
\ns\ds X(t;s,x)=\F(t,s)x+\int_s^t\F(t,\t)b(\t)d\t,\qq
t\in[s,T],\ea\ee
for the solution $X(\cd)\equiv X(\cd\,;s,x,b(\cd))$ of
(\ref{LSDE1}). We have the following result.

\ms

\bf Proposition 2.2. \sl Let {\rm(FD1)} hold. Then the stochastic
fundamental matrix $\F(\cd\,,\cd)$ of $\{A_0(\cd),A_1(\cd)\}$
satisfies the following:
\bel{2.9}\F(t,s)x\ge0,\qq\forall x\ge0,\q 0\le s\le t\le T,~\as,\ee
if and only if
\bel{A0 in R*}A_0(t)\in\dbR^{n\times n}_{*+},\qq t\in[0,T],~~\as,\ee
and
\bel{A1 in Rd}A_1(t)\in\dbR^{n\times n}_d,\qq t\in[0,T],~~\as\ee
Consequently, in this case, for any $(s,x)\in[0,T)\times\dbR^n$ and
$b(\cd)\in L^2_\dbF(\O;L^1(s,T;\dbR^n))$ with
\bel{2.10}x\ge0,\q b(t)\ge0,\qq\ae t\in[s,T],~\as,\ee
the unique solution $X(\cd)\equiv X(\cd\,;s,x,b(\cd))\in
L^2_\dbF(\O;C([s,T];\dbR^n))$ of linear FSDE $(\ref{LSDE1})$
corresponding to $(x,b(\cd))$ on $[s,T]$ satisfies the following:
\bel{ge0}X(t)\ge0,\qq\forall t\in[s,T],\q\as\ee

\ms

\rm

The above result should be known (at least for the case $n=1$). For
reader's convenience, we provide a proof here, which is
straightforward.

\ms

\it Proof. \rm Sufficiency. Let $X(\cd)\equiv X(\cd\,;s,x,0)$ be the
solution to linear FSDE (\ref{LSDE1}) with
$(s,x)\in[0,T)\times\dbR^n$ and $b(\cd)=0$. Then
$$X(t)=\F(t,s)x,\qq0\le s\le t\le T.$$
It suffices to show that $x\le0$ implies
\bel{le0}X(t)\le0,\qq t\in[s,T],\q\as\ee
To prove (\ref{le0}), we define a convex function
$$f(x)=\sum_{i=1}^n(x_i^+)^2,\qq\forall
x=(x_1,x_2,\cds,x_n)\in\dbR^n,$$
where $a^+=\max\{a,0\}$ for any $a\in\dbR$. Applying It\^o's formula
to $f(X(t))$, we get
$$\ba{ll}
\ns\ds f(X(t))-f(x)=\int_s^t\[\lan
f_x(X(\t)),A_0(\t)X(\t)\ran+{1\over2}\lan
f_{xx}(X(\t))A_1(\t)X(\t),A_1(\t)X(\t)\ran\]d\t\\
\ns\ds\qq\qq\qq\qq\qq+\int_s^t\lan f_x(X(\t)),A_1(\t)X(\t)\ran
dW(\t).\ea$$
Let us observe the following: (noting $A_0(\t)\in\dbR^{n\times
n}_{*+}$)
$$\ba{ll}
\ns\ds\lan f_x(X(\t)),A_0(\t)X(\t)\ran=\sum_{i,j=1}^n2X_i(\t)^+\lan
e_i,A_0(\t)e_j\ran X_j(\t)\\
\ns\ds\qq=\sum_{i=1}^n2X_i(\t)^+\lan e_i,A_0(\t)e_i\ran
X_i(\t)+\sum_{i\ne j}2X_i(\t)^+\lan
e_i,A_0(\t)e_j\ran X_j(\t)\\
\ns\ds\qq\le\sum_{i=1}^n2[X_i(\t)^+]^2\lan
e_i,A_0(\t)e_i\ran+\sum_{i\ne j}2\lan e_i,A_0(\t)e_j\ran X_i(\t)^+
X_j(\t)^+\le Kf(X(\t)).\ea$$
Next, we have (noting $A_1(\cd)$ and $f_{xx}(\cd)$ are diagonal)
$$\ba{ll}
\ns\ds{1\over2}\dbE\lan
f_{xx}(X(\t))A_1(\t)X(\t),A_1(\t)X(\t)\ran={1\over2}\dbE\sum_{i=1}^n
I_{(X_i(\t)\ge0)}\(\lan
A_1(\t)e_i,e_i\ran X_i(\t)\)^2\\
\ns\ds={1\over2}\dbE\sum_{i=1}^n\lan
A_1(\t)e_i,e_i\ran{}^2[X_i(\t)^+]^2\le Kf(X(\t)).\ea$$
Consequently,
$$\dbE f(X(t))\le f(x)+K\int_s^t\dbE f(X(\t))d\t,\qq t\in[s,T].$$
Hence, by Gronwall's inequality, we obtain
$$\sum_{i=1}^n\dbE|X_i(t)^+|^2\le K\sum_{i=1}^n|x_i^+|^2,\qq t\in[s,T].$$
Therefore, if $x\le0$, then
$$\sum_{i=1}^n\dbE|X_i(t)^+|^2=0,\qq\forall t\in[s,T].$$
This leads to (\ref{le0}).

\ms

Necessity. Let
$$\O^0_{ij}(s)=\big\{\o\in\O\bigm|\lan A_0(s)e_i,e_j\ran<0\big\},\qq1\le i,j\le n,~s\in[0,T].$$
Suppose (\ref{A0 in R*}) fails. Then for some $i\ne j$, and some
$s\in[0,T)$,
$$\dbP\big(\O^0_{ij}(s)\big)>0,$$
i.e., the $(j,i)$-th (off-diagonal) entry of $A_0(s)$ is not almost
surely nonnegative. Let $X(\cd)\equiv X(\cd\,;s,e_i,0)$ be the
solution to linear FSDE (\ref{LSDE1}) with $(s,x,b(\cd))=(s,e_i,0)$.
Then
$$\ba{ll}
\ns\ds\dbE\big[X_j(t)I_{\O^0_{ij}(s)}\big]=\dbE\big[\lan
X(t),e_j\ran I_{\O^0_{ij}(s)}\big]=\int_s^t\dbE[I_{\O^0_{ij}(s)}\lan
A_0(\t)X(\t), e_j\ran\big]d\t\\
\ns\ds\qq\qq\qq\q=\dbE[\lan A_0(s)e_i,e_j\ran
I_{\O^0_{ij}(s)}\big](t-s)+o(t-s)<0,\ea$$
for $t-s>0$ small. Thus, $X(t)=\F(t,s)e_i\ge0$ fails for some
$t\in[s,T]$ that is close to $s$. This shows that (\ref{A0 in R*})
is necessary.

\ms

Next, suppose (\ref{A1 in Rd}) fails, i.e.,
$$\dbP\(\lan A_1(s)e_i,e_j\ran\ne0\)>0,$$
for some $i\ne j$, and $s\in[0,T)$, i.e., the $(j,i)$-th
(off-diagonal) entry of $A_1(s)$ is not identically equal to zero.
Let $\F_0(\cd\,,\cd)$ be the fundamental matrix of $A_0(\cd)$, i.e.,
$$\F_0(t,s)=I+\int_s^tA_0(\t)\F_0(\t,s)d\t,\qq0\le s\le t\le T.$$
Then $\F_0(\cd,\cd)^{-1}$ satisfies
$$\F_0(t,s)^{-1}=I-\int_s^t\F_0(\t,s)^{-1}A_0(\t)d\t,\qq0\le s\le t\le T.$$
Hence,
$$|\F_0(t,s)^{-1}-I|\le K(t-s),\qq0\le s\le t\le T,\q\as$$
Now, let $X(\cd)=X(\cd\,;s,e_i,0)$. Then
$$X(t)=\F_0(t,s)\[e_i+\int_s^t\F_0(\t,s)^{-1}A_1(\t)X(\t)dW(\t)\],\qq0\le s\le
t\le T.$$
Thus, for $j\ne i$,
$$\lan\F_0(t,s)^{-1}X(t),e_j\ran=\int_s^t\lan\F_0(\t,s)^{-1}A_1(\t)X(\t),e_j\ran dW(\t),\qq t\in[s,T].$$
Consequently,
$$\ba{ll}
\ns\ds X_j(t)=\lan X(t),e_j\ran=\lan[I-\F_0(t,s)^{-1}]X(t),e_j\ran
+\int_s^t\lan\F_0(\t,s)^{-1}A_1(\t)X(\t),e_j\ran dW(\t)\\
\ns\ds\qq~=\lan[I-\F_0(t,s)^{-1}]X(t),e_j\ran+\lan
A_1(s)e_i,e_j\ran\[W(t)-W(s)\]\\
\ns\ds\qq\qq\qq+\int_s^t\lan
A_1(s)e_i-\F_0(\t,s)^{-1}A_1(\t)X(\t),e_j\ran dW(\t).\ea$$
Note that
$$\dbE\,|\lan[I-\F_0(t,s)^{-1}]X(t),e_j\ran|\le\dbE\[|I-\F_0(t,s)|\,|X(t)|
\]\le K(t-s).$$
Also,
$$\ba{ll}
\ns\ds\dbE\,\Big|\1n\int_s^t\lan
A_1(s)e_i-\F_0(\t,s)^{-1}A_1(\t)X(\t),e_j\ran dW(\t)\Big|\\
\ns\ds\le
K\dbE\(\int_s^t|A_1(s)e_i-\F_0(\t,s)^{-1}A_1(\t)X(\t)|^2d\t\)^{1\over2}\2n
=o\big((t-s)^{1\over2}\big).\ea$$
Therefore,
\bel{2.13}\left\{\ba{ll}
\ns\ds\big|\dbE X_j(t)\big|^2=o\big(t-s\big),\\
\ns\ds\dbE|X_j(t)|^2=\dbE|\lan
A_1(s)e_i,e_j\ran|^2(t-s)-o(t-s).\ea\right.\ee
If we let
$$X_j(t)^+=X_j(t)\vee0,\qq X_j(t)^-=[-X_j(t)]\vee0,$$
then
$$X_j(t)=X_j(t)^+-X_j(t)^-,\qq|X_j(t)|=X_j(t)^++X_j(t)^-.$$
Consequently, (\ref{2.13}) can be written as
$$\left\{\ba{ll}
\ns\ds\(\dbE X_j(t)^+-\dbE X_j(t)^-\)^2=o(t-s),\\
\ns\ds\dbE[X_j(t)^+]^2+\dbE[X_j(t)^-]^2=\dbE|\lan
A_1(s)e_i,e_j\ran|^2(t-s)-o(t-s).\ea\right.$$
Hence, it is necessary that
$$\dbE[X_j(t)^+]^2,\q\dbE[X_j(t)^-]^2>0,$$
as long as $t-s>0$ is small, which implies
$$\dbP\(X_j(t)<0\)>0,$$
a contradiction. \endpf

\ms

We point out that in the above, the dimension $n\ge1$; and if $n=1$,
conditions (\ref{A0 in R*})--(\ref{A1 in Rd}) are automatically
true.

\ms

Now, let us look at the following general nonlinear FSDEs, in their
integral form: For $i=0,1$,
\bel{FSDE i}
X^i(t)=x^i+\int_s^tb^i(r,X^i(r))dr+\int_s^t\si(r,X^i(r))dW(r),\qq
t\in[s,T].\ee
Note that unlike the drift $b^i(r,x)$, the diffusion $\si(r,x)$ is
independent of $i=0,1$. We introduce the following assumption.

\ms

{\bf(FD2)} For $i=0,1$, the maps
$b^i,\si:[0,T]\times\dbR^n\times\O\to\dbR^n$ are measurable,
$t\mapsto(b^i(t,x),\si(t,x))$ is $\dbF$-progressively measurable,
$x\mapsto(b^i(t,x),\si(t,x))$ is uniformly Lipschitz, and
$t\mapsto(b^i(t,0),\si(t,0))$ is uniformly bounded.

\ms

It is standard that under (FD2), for any
$(s,x^i)\in[0,T)\times\dbR^n$, (\ref{FSDE i}) admits a unique strong
solution $X^i(\cd)\equiv X^i(\cd\,;s,x^i)$. We have the following
comparison theorem.

\ms

\bf Theorem 2.3. \sl Let {\rm(FD2)} hold. Suppose $\bar
b:[0,T]\times\dbR^n\times\O\to\dbR^n$ is measurable, $t\mapsto\bar
b(t,x)$ is $\dbF$-progressively measurable, $\bar b_x(t,x)$ exists
and is uniformly bounded.

\ms

{\rm(i)} Let
\bel{si_x in Rd}\left\{\ba{ll}
\ns\ds\bar b_x(t,x)\in\dbR^n_{*+},\\
\ns\ds\si_x(t,x)\in\dbR^{n\times
n}_d,\ea\right.\qq\forall(t,x)\in[0,T]\times\dbR^n,~\as\ee
Suppose
\bel{b0<bar b<b1}b^0(t,x)\le\bar b(t,x)\le
b^1(t,x),\qq\forall(t,x)\in[0,T]\times\dbR^n,~\as,\ee
Then for any $(s,x^i)\in[0,T)\times\dbR^n$ with
$$x^0\le x^1,$$
the unique solutions $X^i(\cd)\equiv X^i(\cd\,;s,x^i)$ of
$(\ref{FSDE i})$ satisfy
\bel{X0<X1}X^0(t)\le X^1(t),\qq t\in[s,T],~\as\ee

\ms

{\rm(ii)} Suppose
$$b^0(t,x)=\bar
b(t,x)=b^1(t,x),\qq\forall(t,x)\in[0,T]\times\dbR^n,~\as,$$
and $(t,x)\mapsto(\bar b(t,x),\si(t,x))$ is continuous. Then
$(\ref{si_x in Rd})$ is necessary for the conclusion of {\rm(i)} to
hold.

\rm

\ms

\it Proof. \rm (i) Let $\bar x\in\dbR^n$ with
$$x^0\le\bar x\le x^1.$$
Let $\bar X(\cd)$ be the solution to the following FSDE:
$$\bar X(t)=\bar x+\int_s^t\bar b(r,\bar X(r))dr+\int_s^t\si(r,\bar
X(r))dW(r),\qq t\in[s,T].$$
Then
$$\ba{ll}
\ns\ds\bar X(t)-X^0(t)=\bar x-x^0+\int_s^t\[\bar b(r,X^0(r))-
b^0(r,X^0(r))\]dr\\
\ns\ds\qq\qq\qq\qq+\int_s^t\bar b_x(r)[\bar
X(r)-X^0(r)]dr+\int_s^t\si_x(r)[\bar X(r)-X^0(r)]dW(r),\ea$$
where
$$\ba{ll}
\ns\ds\bar b_x(r)=\int_0^1\bar b_x(r,X^0(r)+\l[\bar
X(r)-X^0(r)])d\l\in\dbR^{n\times n}_{*+},\\
\ns\ds\si_x(r)=\int_0^1\si_x(r,X^0(r)+\l[\bar
X(r)-X^0(r)])d\l\in\dbR^{n\times n}_d.\ea$$
Hence, by Proposition 2.2, we obtain
$$X^0(t)\le\bar X(t),\qq t\in[0,T],~\as$$
Similarly, we are able to show that
$$\bar X(t)\le X^1(t),\qq t\in[0,T],~\as$$
Then (\ref{X0<X1}) follows.

\ms

(ii) For any $x,\wt x\in\dbR^n$ and $\wt x\ge 0,$ $\d\ge0$, let
$X^\d(\cd)$ be the solution to the following:
$$X^\d(t)=x+\d\wt x+\int_s^t\bar
b(r,X^\d(r))dr+\int_s^t\si(r,X^\d(r))dW(r),\qq t\in[s,T],$$
and $\wt X(\cd)$ be the solution to the following:
$$\wt X(t)=\wt x+\int_s^t\bar b_x(r,X^0(r))\wt
X(r)dr+\int_s^t\si_x(r,X^0(r))\wt X(r)dW(r),\qq t\in[s,T].$$
Then it is straightforward that
$$\wt X(t)=\lim_{\d\to0}{X^\d(t)-X^0(t)\over\d},\qq t\in[s,T],~\as$$
Hence, the conclusion of (i) implies that
$$\wt X(t)\ge0,\qq\forall t\in[s,T],~\as$$
Then by Proposition 2.2, we must have
$$\bar b_x(r,X^0(r))\in\dbR^{n\times
n}_{*+},\q\si_x(r,X^0(r))\in\dbR^{n\times n}_d,\qq r\in[s,T],~\as$$
Setting $r=s$, we obtain (\ref{si_x in Rd}). \endpf

\subsection{Comparison of adapted solutions to BSDEs.}

We now look at the following $n$-dimensional linear BSDE:
\bel{BSDE0}\left\{\ba{ll}
\ns\ds dY(t)=\[A(t)Y(t)+B(t)Z(t)-g(t)\]dt+Z(t)dW(t),\qq
t\in[0,\t],\\
\ns\ds Y(\t)=\xi,\ea\right.\ee
where $\xi\in L^2_{\cF_\t}(\O;\dbR^n)$, with $\t$ being an
$\dbF$-stopping time taking values in $(0,T]$. The same as (FD1), we
introduce the following hypothesis.

\ms

{\bf(BD1)} The maps $A(\cd),B(\cd)\in
L^\infty_\dbF(\O;C([0,T];\dbR^{n\times n}))$.

\ms

The following is comparable with Proposition 2.2.

\ms

\bf Proposition 2.4. \sl Let {\rm(BD1)} hold. Then for any
$\dbF$-stopping time $\t$ valued in $(0,T]$, any $g(\cd)\in
L^2_\dbF(0,\t;\dbR^n)$ and $\xi\in L^2_{\cF_\t}(\O;\dbR^n)$ with
\bel{xi,g>0}\xi\ge0,\qq g(t)\ge0,\qq\ae t\in[0,\t],~\as,\ee
the adapted solution $(Y(\cd),Z(\cd))$ to BSDE $(\ref{BSDE0})$
satisfies
$$Y(t)\ge0,\qq t\in[0,\t],~\as,$$
if and only if
\bel{A in R*}-A(t)\in\dbR^{n\times n}_{*+},\q B(t)\in\dbR^{n\times
n}_d,\qq t\in[0,T],~~\as\ee

\ms

\it Proof. \rm Sufficiency. Let $s,\t$ be any $\dbF$-stopping times
such that $0\le s<\t\le T$, almost surely. For any $x\in\dbR^n$, let
$X(\cd)$ be the strong solution to the following FSDE:
\bel{FSDE0}\left\{\ba{ll}
\ns\ds dX(t)=-A(t)^TX(t)dt-B(t)^TX(t)dW(t),\qq t\in[s,\t],\\
\ns\ds X(s)=x.\ea\right.\ee
We claim that the following duality relation holds:
\bel{duality1}\lan x,Y(s)\ran=\dbE_s\[\lan
X(\t),\xi\ran+\int_s^\t\lan X(r),g(r)\ran dr\],\ee
where $\dbE_s[\,\cd\,]=\dbE[\cd\,|\,\cF_s]$. In fact, by It\^o's
formula,
$$\ba{ll}
\ns\ds\dbE_s\[\lan X(\t),\xi\ran-\lan
x,Y(s)\ran\]=\dbE_s\int_s^\t\[-\lan A(r)^TX(r),Y(r)\ran+\lan
X(r),A(r)Y(r)+B(r)Z(r)-g(r)\ran\\
\ns\ds\qq\qq\qq\qq\qq\qq\qq\qq\qq-\lan B(r)^TX(r),Z(r)\ran\]
dr=-\dbE_s\int_s^\t\lan X(r),g(r)\ran dr.\ea$$
Hence, (\ref{duality1}) follows.

\ms

Now, for any $x\in\dbR^n_+$, under our conditions, by Proposition
2.2, the solution $X(\cd)$ of (\ref{FSDE0}) satisfies
$$X(t)\ge0,\qq t\in[s,\t],~\as$$
Hence, by duality relation (\ref{duality1}),
$$\lan x,Y(s)\ran=\dbE_s\[\lan X(\t),\xi\ran+\int_s^\t\lan
X(r),g(r)\ran dr\]\ge0,$$
proving our conclusion.

\ms

Necessity. Suppose (\ref{A in R*}) fails. Then, by Proposition 2.2,
for some $i\ne j$ and $s\in[0,T]$, the solution $X(\cd)\equiv
X(\cd\,;s,e_i)$ of (\ref{FSDE0}) satisfies
$$\dbP\big(X_j(\t)<0\big)>0,$$
for some $\t>s$. For such a $\t$, choosing
$\xi=e_jI_{\{X_j(\t)<0\}}$, and $g(\cd)=0$, we have
$$Y_i(s)=\lan e_i,Y(s)\ran=\dbE_s\[\lan X(\t),e_j\ran I_{\{X_j(\t)<0\}}\]
=\dbE_s\[X_j(\t)I_{\{X_j(\t)<0\}}\]<0,$$
a contradiction. \endpf

\ms

We now look at nonlinear $n$-dimensional BSDEs: For $i=0,1$, and
$\dbF$-stopping time $\t$ valued in $[0,T]$,
\bel{BSDE i}
Y^i(t)=\xi^i+\int_t^Tg^i(s,Y^i(s),Z^i(s))ds-\int_t^TZ^i(s)dW(s),\qq
t\in[0,\t].\ee
Let us introduce the following standard assumption.

\ms

{\bf(BD2)} For $i=0,1$, the map
$g^i:[0,T]\times\dbR^n\times\dbR^n\times\O\to\dbR^n$ is measurable,
$s\mapsto g^i(s,y,z)$ is $\dbF$-progressively measurable,
$(y,z)\mapsto g^i(s,y,z)$ is uniformly Lipschitz, $s\mapsto
g^i(s,0,0)$ is uniformly bounded.

\ms

It is well-known that under (BD2), for any $\xi^i\in
L^p_{\cF_T}(\O;\dbR^n)$ (with $p>1$), BSDE (\ref{BSDE i}) admits a
unique adapted solution $(Y^i(\cd),Z^i(\cd))$. Based on Proposition
2.4, we have the following comparison theorem for nonlinear
$n$-dimensional BSDEs.

\ms

\bf Theorem 2.5. \sl Let {\rm(BD2)} hold. Suppose $\bar
g:[0,T]\times\dbR^n\times\dbR^n\times\O\to\dbR^n$ is measurable,
$s\mapsto\bar g(s,y,z)$ is $\dbF$-progressively measurable, $\bar
g_y(s,y,z)$ and $\bar g_z(s,y,z)$ exist and are uniformly bounded.

\ms

{\rm(i)} Suppose
\bel{bar g}\bar g_y(s,y,z)\in\dbR^{n\times n}_{*+},\q \bar
g_z(s,y,z)\in\dbR^{n\times
n}_d,\qq\forall(s,y,z)\in[0,T]\times\dbR^n\times\dbR^n,~\as,\ee
and
\bel{g0<bar g<g1}g^0(s,y,z)\le\bar g(s,y,z)\le
g^1(s,y,z),\qq\forall(s,y,z)\in[0,T]\times\dbR^n\times\dbR^n,~\as\ee
Then for any $\dbF$-stopping time $\t$ valued in $(0,T]$, and any
$\xi^0,\xi^1\in L^2_{\cF_\t}(\O;\dbR^n)$ with
$$\xi^0\le\xi^1,\qq\as,$$
the corresponding adapted solutions $(Y^i(\cd),Z^i(\cd))$ of BSDEs
$(\ref{BSDE i})$ satisfy
\bel{}Y^0(t)\le Y^1(t),\qq t\in[0,\t],~\as\ee

\ms

{\rm(ii)} Suppose
$$g^0(s,y,z)=\bar
g(s,y,z)=g^1(s,y,z),\qq\forall(s,y,z)\in[0,T]\times\dbR^n\times\dbR^n,~\as,$$
and $(s,y,z)\mapsto\bar g(s,y,z)$ is continuous. Then $(\ref{bar
g})$ is necessary for the conclusion of {\rm(i)} to be true.

\ms

\it Proof. \rm (i) Let $\bar\xi\in L^2_{\cF_\t}(\O;\dbR^n)$ such
that
$$\xi^0\le\bar\xi\le\xi^1,\qq\as$$
Let $(\bar Y(\cd),\bar Z(\cd))$ be the adapted solution to the
following BSDE:
$$\bar Y(t)=\bar\xi+\int_t^\t\bar g(s,\bar Y(s),\bar Z(s))ds-\int_t^\t\bar Z(s)
dW(s),\qq t\in[0,\t].$$
Observe
$$\ba{ll}
\ns\ds \bar Y(t)-Y^0(t)=\bar\xi-\xi^0+\int_t^\t\[\bar
g(s,Y^0(s),Z^0(s))
-g^0(s,Y^0(s),Z^0(s))\]ds\\
\ns\ds\qq\qq\qq\q~+\int_t^\t\[A(s)\(\bar Y(s)-Y^0(s)\)+B(s)\(\bar
Z(s)-Z^0(s)\)\]ds-\int_t^\t\(\bar Z(s)-Z^0(s)\)dW(s),\ea$$
where
$$\ba{ll}
\ns\ds A(s)=\int_0^1\bar g_y\big(s,Y^0(s)+\b[\bar
Y(s)-Y^0(s)],\b[\bar
Z(s)-Z^0(s)]\big)d\b\in\dbR^{n\times n}_{*+},\\
\ns\ds B(s)=\int_0^1\bar g_z\big(s,Y^0(s)+\b[\bar
Y(s)-Y^0(s)],\b[\bar Z(s)-Z^0(s)]\big)d\b\in\dbR^{n\times n}_d.\ea$$
Hence, by Proposition 2.4, we obtain our conclusion.

\ms

(ii) For any given deterministic $\t\in[0,T]$, any $\xi,\wt\xi\in
L^2_{\cF_\t}(\O;\dbR^n)$ and $\wt\xi\ge0,$ $\d\ge0$, let
$(Y^\d(\cd),Z^\d(\cd))$ be the adapted solution to the following
BSDE:
$$Y^\d(t)=\xi+\d\wt\xi+\int_t^\t\bar
g(s,Y^\d(s),Z^\d(s))ds-\int_t^\t Z^\d(s)dW(s),\qq t\in[0,\t].$$
In particular,
\bel{Y0}Y^0(t)=\xi+\int_t^\t\bar g(s,Y^0(s),Z^0(s))ds-\int_t^\t
Z^0(s)dW(s),\qq t\in[0,\t].\ee
If we let $(\wt Y(\cd),\wt Z(\cd))$ be the adapted solution to the
following BSDE:
$$\wt Y(t)=\wt\xi+\int_t^\t\(\bar g_y(s,Y^0(s),Z^0(s))\wt Y(s)
+\bar g_z(s,Y^0(s),Z^0(s))\wt Z(s)\)ds-\int_t^\t\wt Z(s)dW(s), \qq
t\in[0,\t],$$
then it is ready to show that
$$\lim_{\d\to0}{Y^\d(t)-Y^0(t)\over\d}=\wt Y(t),\qq\lim_{\d\to0}{Z^\d(t)-Z^0(t)\over\d}
=\wt Z(t),\qq t\in[0,\t],~\as$$
Hence, conclusion of (i) implies that
$$\wt Y(t)\ge0,\qq\forall t\in[0,\t],~\as$$
Consequently, by Proposition 2.4, we obtain
\bel{2.28}\bar g_y(s,Y^0(s),Z^0(s))\in\dbR^{n\times n}_{*+},\qq\bar
g_z(s,Y^0(s),Z^0(s))\in\dbR^{n\times n}_d,\qq s\in[0,\t],~\as,\ee
for the adapted solution $(Y^0(\cd),Z^0(\cd))$ of BSDE (\ref{Y0})
with any $\xi\in L^2_{\cF_\t}(\O;\dbR^n)$. Now let $\t=T$. For any
$s\in[0,T)$ and $y,z\in\dbR^n$, let
$$\xi=y+z[W(T)-W(s)]-\int_s^T\bar g(r,\bar Y^0(r),z)dr,$$
where $\bar Y^0(\cd)$ is the unique solution of (forward) Volterra
integral equation
$$\bar Y^0(t)=y+z[W(t)-W(s)]-\int_s^t\bar g(r,\bar Y^0(r),z)dr,\qq t\in[s,T].$$
Then it is easy to show that $(\bar Y^0(\cd),z)$ is the unique
adapted solution to the following BSDE
$$Y^0(t)=\xi+\int_t^T\bar g(r,Y^0(r),Z^0(r))dr-\int_t^TZ^0(r)dW(r),\qq t\in[s,T].$$
Clearly, $Y^0(s)=y$, and from (\ref{2.28}), we have
$$\bar g_y(s,y,z)\in\dbR^{n\times n}_{*+},\qq\bar
g_z(s,y,z)\in\dbR^{n\times n}_d,\qq\as$$
Hence, (\ref{bar g}) follows. \endpf

\ms

The above result is a slight extension of a relevant one presented
in \cite{Hu-Peng 2006}, allowing $g^0(\cd)$ and $g^1(\cd)$ to be
different for the sufficient part. Note that as long as the map
$\bar g(\cd)$ exists satisfying (\ref{bar g}) and (\ref{g0<bar
g<g1}), we allow the $j$-th component of $g^i(s,y,z)$ to depend on
$k$-th component of $Z$ with $k\ne j$. For example, suppose $\bar
g(\cd)$ satisfies (\ref{bar g}). Then the comparison theorem holds
for the case, say,
$$g^0(s,y,z)=\bar g(s,y,z)-|z|,\qq g^1(s,y,z)=\bar g(s,y,z)+|z|,\qq(s,y,z)
\in[0,T]\times\dbR^n\times\dbR^n.$$
Finally, we point out that our proof is based on the duality and a
corresponding result for linear FSDEs (Proposition 2.2), which is
different from that found in \cite{Hu-Peng 2006}.

\ms

\subsection{Comparison of solutions to FSVIEs.}

\ms

Let us now turn to FSVIEs. We consider the following linear FSVIE:
\bel{FSVIE1}\ba{ll}
\ns\ds X(t)=\f(t)\1n+\2n\int_0^t\2n\(A_0(t,s)X(s)\1n+\1n
b(s)\)ds\1n+\2n\int_0^t\2n\(A_1(t,s)X(s)\1n+\1n\si(s)\)dW(s),\q
t\in[0,T].\ea\ee
Replacing $\f(\cd)$ by
$$\f(\cd)+\int_0^\cd b(s)ds+\int_0^\cd\si(s)dW(s),$$
we see that without loss of generality, it suffices to consider the
following FSVIE:
\bel{FSVIE2}\ba{ll}
\ns\ds
X(t)=\f(t)+\int_0^tA_0(t,s)X(s)ds+\int_0^tA_1(t,s)X(s)dW(s),\qq
t\in[0,T],\ea\ee
namely, we may assume $b(\cd)=\si(\cd)=0$ in (\ref{FSVIE1}). We now
look at a couple of examples which will help us to exclude some
cases for which the comparison theorem may fail in general.

\ms

\bf Example 2.6. \rm Consider the following one-dimensional
equation:
$$X(t)=1-2e^t\int_0^te^{-s}X(s)ds,\qq t\in[0,T].$$
In this case, we have
$$\f(t)=1,\q A_0(t,s)=-2e^{t-s},\q A_1(t,s)=0,\qq\forall(t,s)\in\D^*.$$
To solve it, let
$$x(t)=\int_0^te^{-s}X(s)ds,\qq t\in[0,T].$$
Then
$$\dot x(t)=e^{-t}X(t)=e^{-t}-2x(t),\qq x(0)=0.$$
Hence,
$$x(t)=\int_0^te^{-2(t-s)}e^{-s}ds=e^{-2t}(e^t-1)=e^{-t}-e^{-2t},\qq t\in[0,T].$$
Therefore, the solution $X(\cd)$ is given by
$$X(t)=1-2e^tx(t)=1-2e^t(e^{-t}-e^{-2t})=-1+2e^{-t},\qq t\in[0,T].$$
Consequently, for $T>\ln2$, we have
$$X(T)=-1+2e^{-T}<0.$$

\ms

This example shows that even for the deterministic case, i.e.,
$A_1(\cd\,,\cd)=0$, the comparison of the solutions may fail. This
is mainly due to the fact that $A_0(\cd\,,\cd)$ is negative and
$t\mapsto A_0(t,s)\equiv-2e^{t-s}$ is decreasing.

\ms

\bf Example 2.7. \rm Consider the following one-dimensional FSVIE:
\bel{sde}X(t)=2T-t+\int_0^tX(s)dW(s),\qq t\in[0,T].\ee
Clearly, (\ref{sde}) is a special case of (\ref{FSVIE1}) with
$$\f(t)=2T-t>0,\q A_0(t,s)=0,\q A_1(t,s)=1.$$
Thus, $\f(\cd)$ is (strictly) positive, and both $A_0(\cd\,,\cd)$
and $A_1(\cd\,,\cd)$ are constants. Note that (\ref{sde}) is
equivalent to the following FSDE:
$$\left\{\ba{ll}
\ns\ds dX(t)=-dt+X(t)dW(t),\qq t\in[0,T],\\
\ns\ds X(0)=2T.\ea\right.$$
Therefore, the solution $X(\cd)$ of the above satisfies the
following:
\bel{2.20}X(t)=e^{-{1\over2}t+
W(t)}\[2T-\int_0^te^{{1\over2}s-W(s)}ds\]\le e^{-{1\over2}t+
W(t)}\[2T-\int_0^te^{-W(s)}ds\],\q t\in[0,T].\ee
By the convexity of $\l\mapsto e^\l$, we have
$${1\over t}\int_0^te^{-W(s)}ds\ge e^{-{1\over
t}\int_0^tW(s)ds}.$$
Thus, for any $t>0$, $X(t)<0$ is implied by
$$e^{-{1\over t}\int_0^tW(s)ds}\ge{2T\over t},$$
which is equivalent to the following:
$$-{1\over t}\int_0^tW(s)ds\ge\log{2T\over t}.$$
Since the left hand side of the above is a normal random variable,
we therefore obtain
\bel{claim2}\dbP(X(t)<0)\ge\dbP\(-{1\over
t}\int_0^tW(s)ds\ge\log{K\over t}\)>0.\ee
This means that the comparison theorem fails for this example.

\ms

From the above, we see that when the diffusion is not identically
zero, nonnegativity of the free term $\f(\cd)$ is not enough to
ensure the nonnegativity of the solution $X(\cd)$ to FSVIE
(\ref{FSVIE2}). The main reason for the comparison fails in this
example is due to the fact that $t\mapsto\f(t)$ is decreasing. Next
example is relevant to a result from \cite{Tudor 1989}, and it is
simpler.

\ms

\bf Example 2.8. \rm Consider
\bel{sde2}X(t)=1+\int_0^t{2T-s\over2T-t} X(s)dW(s),\qq t\in[0,T],\ee
We see that the above is a special case of (\ref{FSVIE2}) with
$$\f(t)=1,\q A_0(t,s)=0,\q A_1(t,s)={2T-s\over2T-t}.$$
The main feature of the above is that the diffusion coefficient
$A_1(t,s)$ depends on $(t,s)$ and the variables $t$ and $s$ cannot
be separated, meaning that $A_1(t,s)$ cannot be written as the
product $A_{11}(t)A_{12}(s)$ of some single variable functions
$A_{11}(\cd)$ and $A_{12}(\cd)$. Clearly, the process $\wt
X(t)\equiv(2T-t)X(t)$ satisfies the following FSVIE:
$$\wt X(t)=2T-t+\int_0^t\wt X(s)dW(s),\qq t\in[0,T],$$
which coincides with (\ref{sde}). Hence, by Example 2.5, although
the free term $\f(t)=1>0$ in (\ref{sde2}), we have
$$\dbP\(X(t)<0\)>0,$$
comparison theorem fails for (\ref{sde2}).

\ms

The above example tells us that if $A_1(t,s)$ is not independent of
$t$, even if the free term $\f(\cd)$ is a constant, comparison
theorem could fail in general. Therefore, if a linear FSVIE is
considered for a general comparison theorem, we had better restrict
ourselves to the following type:
\bel{FSVIE3}X(t)=\f(t)+\int_0^tA_0(t,s)X(s)ds+\int_0^tA_1(s)X(s)dW(s),\qq
t\in[0,T].\ee
To present positive results, we introduce the following assumption.

\ms

{\bf(FV1)} The maps $A_0:\D^*\times\O\to\dbR^{n\times n}$ and
$A_1:[0,T]\times\O\to\dbR^{n\times n}$ are measurable and uniformly
bounded. For any $t\in[0,T]$, $s\mapsto(A_0(t,s),A_1(s))$ is
$\dbF$-progressively measurable on $[0,t]$, and for any $s\in[0,T)$,
the map $t\mapsto A_0(t,s)$ is continuous on $[s,T]$.

\ms

We present the following result.

\ms

\bf Proposition 2.9. \sl Let {\rm(FV1)} hold.

\ms

{\rm(i)} Suppose
\bel{2.21}A_0(t,s)\in\dbR^{n\times n}_+,\q
A_1(s)=0,\q\ae(t,s)\in\D^*,~\as\ee
Then for any $\f(\cd)\in L^2_\dbF(0,T;\dbR^n)$ $(\ref{FSVIE3})$
admits a unique solution $X(\cd)\in L^2_\dbF(0,T;\dbR^n)$ and it
satisfies
\bel{5.10}X(t)\ge\f(t)\ge0,\qq t\in[0,T].\ee

{\rm(ii)} Suppose
\bel{2.26}\ba{ll}
\ns\ds A_0(t,s)\in\dbR^{n\times n}_{*+},\q A_1(s)\in\dbR^{n\times
n}_d,\qq\ae(t,s)\in\D^*,~\as\ea\ee
Moreover, there exists a continuous nondecreasing function
$\rho:[0,T]\to[0,\infty)$ with $\rho(0)=0$ such that
\bel{rho}|A_0(t,s)-A_0(t',s)|\le\rho(|t-t'|),\qq
t,t'\in[0,T],~s\in[0,t\land t'],~\as,\ee
and
\bel{A0>A0}A_0(\t,s)-A_0(t,s)\in\dbR^{n\times n}_+,\qq\forall\,0\le
s\le t\le\t\le T,~\as\ee
Then for any $\f(\cd)\in C_\dbF([0,T];L^2(\Omega,\dbR^n))$, with
\bel{f>f}\f(\t)\ge\f(t)\ge0,\qq\forall\,0\le s\le t\le\t\le
T,~\as,\ee
$(\ref{FSVIE3})$ admits a unique solution $X(\cd)\in
C_\dbF([0,T];L^2(\O;\dbR^n))$ and it satisfies:
\bel{X>0}X(t)\ge0,\qq t\in[0,T],~\as\ee

\rm

Note that between the above (i) and (ii), none of them includes the
other. Condition (\ref{2.21}) implies that the map $y\mapsto
A_0(t,s)y$ is nondecreasing (for $y\ge0$); whereas, condition
(\ref{A0>A0}) implies that the map $t\mapsto A_0(t,s)y$ is
nondecreasing. The monotonicity of $\f(\cd)$ is assumed in (ii),
which is not needed in (i). We will encounter a similar situation
for BSVIEs a little later. Also, because of Example 2.6,
$\dbR^{n\times n}_+$ in (\ref{2.21}) cannot be replaced by
$\dbR^{n\times n}_{*+}$.

\ms

\it Proof. \rm (i) Define
$$(\cA X)(t)=\int_0^tA_0(t,s)X(s)ds,\qq t\in[0,T].$$
By our condition, making use of Proposition 2.1, we see that
$$(\cA X)(\cd)\ge0,\qq\forall X(\cd)\in L^2_\dbF(0,T;\dbR^n),~X(\cd)\ge0.$$
Now, we define the following Picard iteration sequence
$$X^0(\cd)=\f(\cd),\qq X^k(\cd)=\f(\cd)+(\cA X^{k-1})(\cd),\qq k\ge1.$$
By induction, it is easy to see that
$$X^k(\cd)\ge\f(\cd),\qq\forall k\ge0.$$
Further,
$$\lim_{k\to\infty}\|X^k(\cd)-X(\cd)\|_{L^2_\dbF(0,T;\dbR^n)}=0,$$
with $X(\cd)$ being the solution to (\ref{FSVIE3}). Then it is easy
to see that (\ref{5.10}) holds.

\ms

(ii) Let $\Pi=\{\t_k,0\le k\le N\}$ be an arbitrary set of finitely
many $\dbF$-stopping times with $0=\t_0<\t_1<\cds<\t_N=T$, and we
define its mesh size by
$$\|\Pi\|=\esssup_{\o\in\O}\max_{1\le k\le N}|\t_k-\t_{k-1}|.$$
Let
$$
A^\Pi_0(t,s)=\sum_{k=0}^{N-1}A_0(\t_k,s)I_{[\t_k,\t_{k+1})}(t),\qq
\f^\Pi(t)=\sum_{k=0}^{N-1} \f(\t_k)I_{[\t_k,\t_{k+1})}(t).$$
Clearly, each $A_0(\t_k,\cd)$ is an $\dbF$-adapted bounded process,
and each $\f(\t_k)$ is an $\cF_{\t_k}$-measurable random variable.
Moreover, for each $k\ge0$,
\bel{}A_0(\t_k,s)\in\dbR^{n\times n}_{*+},\qq
s\in[\t_k,\t_{k+1}),~~\as,\ee
and
\bel{5.13}0\le\f(\t_k)\le\f(\t_{k+1}),\qq\as\ee
Further,
$$\ba{ll}
\ns\ds|A_0^\Pi(t,s)-A_0(t,s)|=\sum_{k=0}^{N-1}|A_0(\t_k,s)-A_0(t,s)|I_{[\t_k,\t_{k+1})}(t)\le\sum_{k=0}^{N-1}\rho(t-\t_k)I_{[\t_k,\t_{k+1})}(t)
\le\rho\big(\|\Pi\|).\ea$$
Now, we let
$X^\Pi(\cd)$ be the solution to the following FSVIE:
\bel{}\ba{ll}
\ns\ds
X^\Pi(t)=\f^\Pi(t)+\int_0^tA_0^\Pi(t,s)X^\Pi(s)ds+\int_0^tA_1(s)X^\Pi(s)dW(s),\qq
t\in[0,T].\ea\ee
Then we can show that
%
%
%
%
%
%
%
%
%
%
%
%
%
%
%
%
%
%
\bel{lim}\lim_{\|\Pi\|\to0}\dbE\[\sup_{t\in[0,T]}|X^\Pi(t)-X(t)|^2\]=0.\ee
We now want to show that
\bel{X Pi>0}X^\Pi(t)\ge0,\qq t\in[0,T],~\as,\ee
which, together with (\ref{lim}) will lead to (\ref{X>0}). To show
(\ref{X Pi>0}), we look at $X^\Pi(\cd)$ on each interval
$[\t_k,\t_{k+1})$, $k=0,1,\cds,N-1$. First, on interval $[0,\t_1)$,
we have
$$
X^\Pi(t)=\f(0)+\int_0^tA_0(0,s)X^\Pi(s)ds+\int_0^tA_1(s)X^\Pi(s)dW(s),$$
which is an FSDE, and $X^\Pi(\cd)$ has continuous paths (on
$[0,\t_1)$). From Proposition 2.2, we have
$$X^\Pi(t)\ge0,\qq t\in[0,\t_1),~\as$$
In particular,
\bel{2.16}\ba{ll}
\ns\ds
X^\Pi(\t_1-0)=\f(0)+\int_0^{\t_1}A_0(0,s)X^\Pi(s)ds+\int_0^{\t_1}A_1(s)X^\Pi(s)dW(s)\ge0.\ea\ee
Next, on $[\t_1,\t_2)$, we have (making use of (\ref{2.16}))
$$\ba{ll}
\ns\ds
X^\Pi(t)=\f(\t_1)+\int_0^{\t_1}A_0(\t_1,s)X^\Pi(s)ds+\int_0^{\t_1}
A_1(s)X^\Pi(s)dW(s)\\
\ns\ds\qq\qq+\int_{\t_1}^tA_0(\t_1,s)X^\Pi(s)ds
+\int_{\t_1}^tA_1(s)X^\Pi(s)dW(s)\\
\ns\ds\qq=\f(\t_1)-\f(0)+X^\Pi(\t_1-0)+\int_0^{\t_1}\(A_0(\t_1,s)-A_0(0,s)\)X^\Pi(s)ds\\
\ns\ds\qq\qq+\int_{\t_1}^t\(A_0(\t_1,s)X^\Pi(s)\)ds
+\int_{\t_1}^tA_1(s)X^\Pi(s)dW(s)\\
\ns\ds\qq\equiv\wt X(\t_1)+\int_{\t_1}^tA_0(\t_1,s)X^\Pi(s)ds
+\int_{\t_1}^tA_1(s)X^\Pi(s)dW(s),\ea$$
where, by (\ref{f>f}) and (\ref{2.16}),
$$\ba{ll}
\ns\ds\wt
X(\t_1)\equiv\f(\t_1)-\f(0)+X^\Pi(\t_1-0)+\int_0^{\t_1}\(A_0(\t_1,s)-A_0(0,s)\)X^\Pi(s)ds\ge0.\ea$$
Hence, one obtains
$$X^\Pi(t)\ge0,\qq t\in[\t_1,\t_2).$$
By induction, we obtain (\ref{X Pi>0}). \endpf

\ms

\ms

Based on the above result, it is not very hard for us to present
comparison theorems for nonlinear FSVIEs. We prefer not to give the
details here. One can cook up that by following the relevant details
for BSVIEs which will be presented in the following section. To
conclude this section, we present an example showing that in the
case $A_1(\cd)\ne0$, as long as $t\mapsto A_0(t,\cd)$ is not
nondecreasing in the sense of (\ref{A0>A0}), even if
$A_0(t,s)\in\dbR^{n\times n}_+$, comparison theorem might still fail
as well.

\ms

\bf Example 2.10. \rm Consider the following FSVIE:
$$X(t)=1+\int_0^tI_{[0,\t]}(t)X(s)ds+\int_0^tX(s)dW(s),\qq
t\in[0,T],$$
where $\t\in(0,T)$. Clearly, the above is a special case of
(\ref{FSVIE2}) with
$$\f(t)=1,\q A_0(t,s)=I_{[0,\t]}(t),\q A_1(s)=1.$$
Thus, $t\mapsto A_0(t,s)$ is not nondecreasing. Let us solve this
FSVIE. On $[0,\t)$, we have
$$X(t)=1+\int_0^tX(s)ds+\int_0^tX(s)dW(s),$$
which is equivalent to the following:
$$dX(t)=X(t)dt+X(t)dW(t),\qq X(0)=1.$$
Hence,
$$X(t)=e^{{t\over2}+W(t)},\qq t\in[0,\t).$$
On $[\t,T]$, we have
$$\ba{ll}
\ns\ds
X(t)=1+\int_0^tX(s)dW(s)=1+\int_0^\t X(s)dW(s)+\int_\t^tX(s)dW(s)\\
\ns\ds\qq=X(\t-0)-\int_0^\t X(s)ds+\int_\t^tX(s)dW(s),\ea$$
which is equivalent to the following:
$$dX(t)=X(t)dW(t),\qq X(\t+0)=X(\t-0)-\int_0^\t X(s)ds.$$
Hence,
$$\ba{ll}
\ns\ds X(t)=e^{-{t-\t\over2}+W(t)-W(\t)}\[X(\t-0)-\int_0^\t X(s)ds\]\\
\ns\ds\qq=e^{-{t-\t\over2}+W(t)-W(\t)}\[e^{{\t\over2}+W(\t)}
-\int_0^\t e^{{s\over2}+W(s)}ds\].\ea$$
Then $X(t)<0$ for $t\in[\t,T]$ if and only if
$$e^{{\t\over2}+W(\t)}
<\int_0^\t e^{{s\over2}+W(s)}ds.$$
By convexity of $\l\mapsto e^\l$, we have
$${1\over\t}\int_0^\t e^{{s\over2}+W(s)}ds\ge
e^{{1\over\t}\int_0^\t[{s\over2}+W(s)]ds}.$$
Hence, $X(t)<0$ for some $t\in[\t,T]$ is implied by
$$e^{{\t\over2}+W(\t)}<\t e^{{1\over\t}\int_0^\t({s\over2}+W(s))ds},$$
which is equivalent to
$$\ba{ll}
\ns\ds{\t\over2}+W(\t)<\ln\t+{1\over\t}\int_0^\t\({s\over2}+W(s)
\)ds=\ln\t+{\t\over4}+{1\over\t}\int_0^\t W(s)ds\\
\ns\ds\qq\qq\q=\ln\t+{\t\over4}+{1\over\t}sW(s)\Big|_0^\t -\int_0^\t
sdW(s)=\ln\t+{\t\over4}+W(\t)-\int_0^\t sdW(s).\ea$$
This is further equivalent to the following:
$$\int_0^\t sdW(s)<\ln\t-{\t\over4}.$$
The left hand side of the above is a normal random variable. Hence,
$$\dbP\(\int_0^\t sdW(s)<\ln\t-{\t\over4}\)>0,$$
which implies
$$\dbP\(X(t)<0\)>0,\qq t\in(\t,T].$$
Although in the above, $A_0(\cd\,,\cd)$ is discontinuous, it is not
hard for us to replace it by a continuous one and still have the
same conclusion.

\ms

\section{Comparison Theorems for BSVIEs.}

In this section we consider various comparison theorems for BSVIEs.

\ms

\subsection{Comparison for adapted solutions.}

We first consider the following type BSVIEs: For $i=0,1$,
\bel{BSVIE5}Y^i(t)=\psi^i(t)+\int_t^Tg^i(t,s,Y^i(s),Z^i(t,s))ds
-\int_t^TZ^i(t,s)dW(s),\q t\in[0,T].\ee
The key feature here is that the generator $g^i(\cd)$ is independent
of $Z^i(s,t)$. For any adapted solution
$(Y^i(\cd),Z^i(\cd\,,\cd))\in\cH^p_\D[0,T]$ of the above, we need
only the values $Z^i(t,s)$ of $Z^i(\cd\,,\cd)$ for $0\le t\le s\le
T$, and the values $Z^i(t,s)$ of $Z^i(\cd\,,\cd)$ for $0\le s<t\le
T$ are irrelevant. Consequently, the notion of M-solution is not
necessary for BSVIE of form (\ref{BSVIE5}). For the generator
$g(\cd)$ of BSVIE (\ref{BSVIE5}), we adopt the following assumption.

\ms

{\bf(BV1)} Let $g^i:\D\times\dbR^n\times\dbR^n\times\O\to\dbR^n$ be
measurable such that $s\mapsto g^i(t,s,y,z)$ is $\dbF$-progressively
measurable, $(y,z)\mapsto g^i(t,s,y,z)$ is uniformly Lipschitz,
$(t,s)\mapsto g^i(t,s,0,0)$ is uniformly bounded.

\ms

It is known that under (BV1), for any $\psi^i(\cd)\in
C_{\cF_T}([0,T];L^2(\O;\dbR^n))$, BSVIE (\ref{BSVIE5}) admits a
unique adapted solution $(Y^i(\cd),Z^i(\cd\,,\cd))\in\cH^2_\D[0,T]$.
We want to look at if a proper comparison between $Y^1(\cd)$ and
$Y^0(\cd)$ holds under certain additional conditions on $g^i(\cd)$
and $\psi^i(\cd)$. To begin with, let us first look at the following
simple BSVIEs: For $i=0,1$,
\bel{BSVIE6_1}Y^i(t)=\psi^i(t)+\int_t^Tg^i(t,s,Z^i(t,s))ds-\int_t^TZ^i(t,s)dW(s),\qq
t\in[0,T],\ee
with the generators $g^i(\cd)$ are independent of $Y^i(s)$. We have
the following result.

\ms

\bf Proposition 3.1. \sl For $i=0,1$, let
$g^i:\D\times\dbR^n\times\O\to\dbR^n$ satisfy {\rm(BV1)}. Moreover,
\bel{}g^0(t,s,z)\le
g^1(t,s,z),\qq\forall(t,s,z)\in\D\times\dbR^n,~\as\ee
and for either $i=0$ or $i=1$, $g^i_z(t,s,z)$ exists and
\bel{}g^i_z(t,s,z)\in\dbR^{n\times
n}_d,\qq(t,s,z)\in\D\times\dbR^n,~\as\ee
Then the adapted solutions
$(Y^i(\cd),Z^i(\cd\,,\cd))\in\cH^2_\D[0,T]$ of BSVIE
$(\ref{BSVIE6_1})$ with
\bel{}\psi^0(t)\le\psi^1(t),\qq t\in[0,T],~\as,\ee
satisfies
\bel{Y>Y3.1_1}Y^0(t)\le Y^1(t),\qq t\in [0,T],~\as\ee

\ms

\it Proof. \rm Fixed $t\in [0,T]$. For $i=0,1$, let
$(\l^i(t,\cd),\m^i(t,\cd))$ be the adapted solution to the following
BSDE:
$$\lambda ^i(t,r)=\psi^i(t)+\int_r^Tg^i(t,s,\mu
^i(t,s))ds-\int_r^T\mu ^i(t,s)dW(s),\qq r\in[t,T].$$
By Theorem 2.5, we have that
\bel{l>l}\l^0(t,r)\le\l^1(t,r),\qq r\in[t,T],~\as\ee
By setting
\bel{YZ}Y^i(t)=\l^i(t,t),\qq
Z^i(t,s)=\m^i(t,s),\qq\forall(t,s)\in\D,\ee
we see that $(Y^i(\cd),Z^i(\cd\,,\cd))$ is the adapted solution to
the BSVIE (\ref{BSVIE6_1}). Then (\ref{Y>Y3.1_1}) follows from
(\ref{l>l}). \endpf

\ms

Returning to BSVIEs (\ref{BSVIE5}), we have the following result.

\ms

\bf Theorem 3.2. \sl Let {\rm(BV1)} hold. Suppose $\bar
g:\D\times\dbR^\times\dbR^n\times\O\to\dbR^n$ is measurable,
$s\mapsto\bar g(t,s,y,z)$ is $\dbF$-progressively measurable,
$(y,z)\mapsto\bar g(t,s,y,z)$ is uniformly Lipschitz, $y\mapsto\bar
g(t,s,y,z)$ is nondecreasing, such that
\bel{g0<bar g<g1}g^0(t,s,y,z)\le\bar g(t,s,y,z)\le
g^1(t,s,y,z),\qq(t,s,y,z)\in\D\times\dbR^n\times\dbR^n,~\as\ee
Moreover, $\bar g_z(t,s,y,z)$ exists and
\bel{bar gz in Rd}\bar g_z(t,s,y,z)\in\dbR^{n\times
n}_d,\qq(t,s,y,z)\in\D\times\dbR^n\times\dbR^n,~\as\ee
Then for any $\psi^i(\cd)\in C_{\cF_T}([0,T];L^2(\O;\dbR^n))$
satisfying
\bel{psi>psi}\psi^0(t)\le\psi^1(t),\qq t\in[0,T],~\as,\ee
the corresponding unique adapted solution
$(Y^i(\cd),Z^i(\cd\,,\cd))\in\cH^2_\D[0,T]$ of BSVIE
$(\ref{BSVIE5})$ satisfy
\bel{Y>Y3.1}Y^0(t)\le Y^1(t),\qq t\in [0,T],~\as\ee

\ms

\it Proof. \rm Let $\bar\psi(\cd)\in
C_{\cF_T}([0,T];L^2(\O;\dbR^n))$ such that
$$\psi^0(t)\le\bar\psi(t)\le\psi^1(t),\qq t\in[0,T],~\as$$
Let $(\bar Y(\cd),\bar Z(\cd\,,\cd))$ be the adapted solution to the
following:
$$\bar Y(t)=\bar\psi(t)+\int_t^T\bar g(t,s,\bar Y(s),\bar
Z(t,s))ds-\int_t^T\bar Z(t,s)dW(s),\qq t\in[0,T].$$
Set $\wt Y_0(\cd)=Y^0(\cd)$ and consider the following BSVIE:
$$\wt Y_1(t)=\bar\psi(t)+\int_t^T\bar g(t,s,\wt Y_0(s),\wt Z_1(t,s))ds
-\int_t^T\wt Z_1(t,s)dW(s),\q t\in[0,T].$$
Let $(\wt Y_1(\cd),\wt Z_1(\cd\,,\cd))\in\cH^2_\D[0,T]$ be the
unique adapted solution to the above. Since
$$\left\{\ba{ll}
\ns\ds\bar g(t,s,\wt Y_0(s),z)\le
g^1(t,s,\wt Y_0(s),z),\qq(t,s,z)\in\D\times\dbR^n,~\as,\\
\ns\ds\bar g_z(t,s,\wt Y_0(s),z)\in\dbR^{n\times
n}_d,\qq\qq(t,s,z)\in\D\times\dbR^n,~\as,\\
\ns\ds\bar\psi(t)\le\psi^1(t),\qq t\in[0,T],~\as\ea\right.$$
By Proposition 3.1, we obtain that
$$\wt Y_1(t)\le\wt Y_0(t),\qq t\in[0,T].$$
Next, we consider the following BSVIE:
$$\wt Y_2(t)=\bar\psi(t)+\int_t^T\bar g(t,s,\wt Y_1(s),\wt Z_2(t,s))ds
-\int_t^T\wt Z_2(t,s)dW(s),\qq t\in[0,T],$$
and let $(\wt Y_2(\cd),\wt Z_2(\cd\,,\cd))\in\cH^2_\D[0,T]$ be the
adapted solution to the above. Now, since $y\mapsto\bar g(t,s,y,z)$
is nondecreasing, we have
$$\bar g(t,s,\wt Y_1(s),z)\le\bar g(t,s,\wt Y_0(s),z),\qq\forall(t,s,z)\in\D\times\dbR^n.$$
Hence, similar to the above, we obtain
$$\wt Y_2(t)\le\wt Y_1(t),\qq t\in[0,T],~\as$$
By induction, we can construct a sequence $\{(\wt Y_k(\cd),\wt
Z_k(\cd\,,\cd))\}_{k\ge1}\subseteq\cH^2_\D[0,T]$ such that
$$\wt Y_k(t)=\bar\psi(t)+\int_t^T\bar g(t,s,\wt Y_{k-1}(s),\wt Z_k(t,s))ds
-\int_t^T\wt Z_k(t,s)dW(s),\qq t\in[0,T],$$
and
\bel{Y0>Y1}Y^1(t)=\wt Y_0(t)\ge\wt Y_1(t)\ge\wt Y_2(t)\cds,\qq
t\in[0,T],~\as\ee
Next we will show that the sequence $\{(\wt Y_k(\cd),\wt
Z_k(\cd\,,\cd))\}_{k\ge1}$ is Cauchy in $\cH^2_\D[0,T]$. To show
this, we introduce an equivalent norm of $\cH^2_\D[0,T]$ as
$$\|(y(\cd),z(\cd\,,\cd))\|_{\cH^2_{\Delta}[0,T]}^2=\dbE\int_0^Te^{\b t}|y(t)|^2dt
+\dbE\int_0^Te^{\b t}\int_t^T|z(t,s)|^2dsdt,$$
with $(y(\cd),z(\cd\,,\cd))\in\cH^2_\D[0,T]$, and $\b$ being a
constant undetermined. By utilizing a stability estimate in
\cite{Yong 2008}, we have
\bel{}\ba{ll}
\ns\ds\dbE|\wt Y_k(t)-\wt Y_\ell(t)|^2+\dbE\int_t^T|\wt Z_k(t,s)
-\wt Z_\ell(t,s)|^2ds\\
\ns\ds\le K\dbE\(\int_t^T|\bar g(t,s,\wt Y_{k-1}(s),\wt
Z_k(t,s))-\bar g(t,s,\wt Y_{\ell-1}(s),\wt Z_k(t,s))|ds\)^2.\ea\ee
Consequently, we arrive at
\bel{}\ba{ll}
\ns\ds\dbE\int_0^Te^{\b t}|\wt Y_k(t)-\wt
Y_\ell(t)|^2dt+\dbE\int_0^Te^{\b
t}\(\int_t^T|\wt Z_k(t,s)-\wt Z_\ell(t,s)|^2ds\)dt\\
\ns\ds\le K\dbE\int_0^Te^{\b t}\(\int_t^T|\bar g(t,s,\wt
Y_{k-1}(s),\wt Z_k(t,s))-\bar g(t,s,\wt Y_{\ell-1}(s),\wt Z_k(t,s))|ds\)^2dt\\
\ns\ds\le K\dbE\int_0^Te^{\b t}\(\int_t^T|\wt Y_{k-1}(s)-\wt
Y_{\ell-1}(s)|
ds\)^2dt\\
\ns\ds\le K\dbE\int_0^T|\wt Y_{k-1}(s)-\wt
Y_{\ell-1}(s)|^2ds\int_0^s e^{\b t}dt\le{K\over\b}\dbE\int_0^Te^{\b
s}|\wt Y_{k-1}(s)-\wt Y_{\ell-1}(s)|^2ds.\ea\ee
Note that the constant $K>0$ in the above can be chosen independent
of $\b>0$. Thus by choosing a $\b$ such that ${K\over\b}<1$, we
obtain immediately that $\{(\wt Y_k(\cd),\wt
X_k(\cd\,,\cd))\}_{k\ge1}$ is Cauchy in $\cH^2_\D[0,T]$. Hence,
there exists a $(\wt Y(\cd),\wt Z(\cd\,,\cd))\in\cH^2_\D[0,T]$ such
that
$$\lim_{k\to\infty}\[
\dbE\int_0^T|\wt Y^k(t)-\wt Y(t)|^2dt+\dbE\int_0^Te^{\b
t}\(\int_t^T|\wt Z^k(t,s)-\wt Z(t,s)|^2ds\)dt\]=0,$$
and
$$\wt Y(t)=\bar\psi(t)+\int_t^T\bar g(t,s,\wt Y(s),\wt Z(t,s))ds
-\int_t^T\wt Z(t,s)dW(s),\qq t\in[0,T].$$
By uniqueness, we have
$$\bar Y(t)=\wt Y(t)\le\wt Y_0(t)=Y^1(t),\qq t\in[0,T],~\as$$
Similarly, we can prove that
$$Y^0(t)\le\bar Y(t),\qq t\in[0,T],~\as$$
Therefore, our conclusion follows. \endpf

\ms

It is easy to cook up an example for which $y\mapsto g^i(t,s,y,z)$
is not nondecreasing for $i=0,1$, but a $\bar g(\cd)$ satisfying
conditions of Theorem 3.2 can be constructed. For example,
$$g^0(t,s,y,z)\equiv\sin y\le 1\equiv\bar g(t,s,y,z)\le 2+\cos y\equiv
g^1(t,s,y,z).$$
Condition (\ref{g0<bar g<g1}) means that in the tube
$$\Big\{\[g^0(t,s,y,z),g^1(t,s,y,z)\]\bigm|(t,s,y,z)\in\D\times\dbR^n\times
\dbR^n\Big\},$$
there exists a selection $\bar g(t,s,y,z)$ which is nondecreasing in
$y$, and (\ref{bar gz in Rd}) is satisfied. Therefore, the condition
assumed in Theorem 3.2 is a kind of {\it generalized nondecreasing
condition} for the maps $y\mapsto g^i(t,s,y,z)$, although these maps
themselves are not necessarily nondecreasing. Consequently, it is
expected that condition (\ref{g0<bar g<g1}) excludes many other
situations. To see that, let us look at two examples.

\ms

\bf Example 3.3. \rm Consider one-dimensional linear BSVIE
$$Y(t)=t-\int_t^TY(s)ds-\int_t^TZ(t,s)dW(s),\qq t\in[0,T].$$
It is clear that if $(Y(\cd),Z(\cd\,,\cd))\in\cH^2_\D[0,T]$ is the
adapted solution, then $Z(\cd\,,\cd)=0$ and
$$Y(t)=e^{t-T}(T+1)-1,\qq t\in[0,T].$$
Consequently,
$$Y(t)<0,\qq t\in[0,T-\ln(T+1)].$$
Therefore, comparison theorem fails for this example. This example
corresponds to the case
$$g^i(t,s,y,z)=-y,\q i=0,1,\qq\psi^1(t)=t,\q\psi^0(t)=0.$$
Apparently, $\bar g(\cd)$ satisfying the conditions in Theorem 3.2
does not exist.

\ms

\bf Example 3.4. \rm Consider
$$Y(t)=1+\int_t^T(t-1)Y(s)ds-\int_t^TZ(t,s)dW(s),\qq t\in[0,T].$$
Again, if $(Y(\cd),Z(\cd\,,\cd))\in\cH^2_\D[0,T]$ the adapted
solution, then $Z(\cd\,,\cd)=0$. Now, we denote
$$y(t)=\int_t^TY(s)ds.$$
Then
$$\dot y(t)=-Y(t)=-1-(t-1)\int_t^TY(s)ds=-1-(t-1)y(t).$$
Hence,
$$\ba{ll}
\ns\ds0=y(T)=e^{-\int_t^T(s-1)ds}y(t)-\int_t^Te^{-\int_\t^T(s-1)ds}d\t.\ea$$
This yields
$$y(t)=\int_t^Te^{\int_t^\t(s-1)ds}d\t=\int_t^Te^{{1\over2}[(\t-1)^2-(t-1)^2]}d\t.$$
Therefore,
$$Y(t)=1+(t-1)y(t)=1+(t-1)\int_t^Te^{{1\over2}[\t^2-t^2-2(\t-t)]}d\t.$$
Consequently,
$$Y(0)=1-\int_0^Te^{{1\over2}\t^2-\t}d\t<0,$$
provided $T>0$ is large.
Thus, comparison theorem fails for this example as well. This
example corresponds to the case
$$g^i(t,s,y,z)=(t-1)y,\q i=0,1,\qq\psi^1(t)=1,\q\psi^0(t)=0.$$
Again, for this example, the generator $\bar g(\cd)$ satisfying the
conditions in Theorem 3.2 does not exist.

\ms

Let us take a closer look at the above two examples. In Example 3.3,
$t\mapsto\psi^1(t)$ is increasing, and in Example 3.4, $t\mapsto
g^i(t,s,y,z)$ is increasing for $y>0$. In a certain sense, these
conditions actually prevent the comparison theorem from being true
for these examples. On the other hand, we keep in mind that when
$\psi(t)$ and $g^i(t,s,y,z)$ are independent of $t$, the above two
situations do not appear. Hence, it is natural to ask if comparison
theorem remains when $\psi(t)$ and $g^i(t,s,y,z)$ do depend on $t$,
and the generalized nondecreasing condition (\ref{g0<bar g<g1}) is
not assumed. The answer is positive. Before we state and prove a
general positive result, let us look at the following example.

\ms

\bf Example 3.5. \rm Consider the following BSVIE:
$$Y(t)=\int_t^T\[s-t-Y(s)\]ds-\int_t^TZ(t,s)dW(s),\qq t\in[0,T].$$
In this case, we have
$$\psi(t)\equiv0,\quad g(t,s,y,z)=s-t-y.$$
Thus, condition of Theorem 3.2 fails. However, it is easy to check
that the unique adapted solution $(Y(\cd),Z(\cd\,,\cd))$ is given by
$$Y(s)=e^{s-T}+T-s-1,\q Z(t,s)=0.$$
is the unique solution here. Clearly,
$$Y(s)\ge0,\qq s\in[0,T],$$
comparison theorem holds. Note that in this case, $t\mapsto
g(t,s,y,z)$ is nondecreasing. On the other hand, the BSVIE is
equivalent to the following:
$$Y(t)={(T-t)^2\over2}-\int_t^TY(s)ds-\int_t^TZ(t,s)dW(s),\qq
t\in[0,T],$$
with
$$\psi(t)={(T-t)^2\over2},\qq g(t,s,y)=-y.$$
For this, we have that $t\mapsto\psi(t)$ is non-increasing.

\ms

Inspired by the above example, we see that without condition
(\ref{g0<bar g<g1}), one might still have comparison theorem. We now
establish such kind of results. Let us begin with a result for
linear BSVIEs. More precisely, we consider the following linear
BSVIE:
\bel{3.14}Y(t)=\psi(t)+\int_t^T\[A(t,s)Y(s)+B(s)Z(t,s)\]ds-\int_t^TZ(t,s)dW(s),\q
t\in[0,T].\ee
Note that the coefficient $B(s)$ of $Z(t,s)$ is independent of $t$.
We have the following theorem.

\ms

\bf Theorem 3.6. \sl Let $A:\D\times\O\to\dbR^{n\times n}$ and
$B:[0,T]\times\O\to\dbR^{n\times n}$ be uniformly bounded, with
$B(\cd)$ being $\dbF$-progressively measurable, for each
$t\in[0,T]$, $s\mapsto A(t,s)$ being $\dbF$-progressively
measurable, and for each $s\in[0,T]$, $t\mapsto A(t,s)$ being
continuous. Moreover,
\bel{3.18}A(t,s)\in\dbR^{n\times n}_{*+},\qq(t,s)\in\D,~\as,\ee
\bel{A-A>0}A(t,s)-A(\t,s)\in\dbR^{n\times n}_+,\qq0\le t\le\t\le
s\le T,~\as,\ee
\bel{}B(s)\in\dbR^{n\times n}_d,\qq s\in[0,T],~\as\ee
Then for any $\psi(\cd)\in C_{\cF_T}([0,T];L^2(\O;\dbR^n))$ with
\bel{}\psi(t)\ge\psi(s)\ge0,\qq0\le t\le s\le T,~\as,\ee
the adapted solution $Y(\cd),Z(\cd\,,\cd))$ of linear BSVIE
$(\ref{3.14})$ satisfies the following:
\bel{3.22}Y(t)\ge0,\qq t\in[0,T],~\as\ee

\ms

\rm

We point out that $A(t,s)$ satisfying (\ref{3.18}) (which is always
true if $n=1$) is not necessarily in $\dbR^{n\times n}_+$.
Therefore, the map $y\mapsto A(t,s)y$ is not necessarily
nondecreasing. Also, when $A(t,s)$ is independent of $t$,
(\ref{A-A>0}) is automatically true.

\ms

\it Proof. \rm Let
$$A(t,s)=\sum_{k=1}^NA_k(s)I_{(t_{k-1},t_k]}(t),\q
\psi(t)=\sum_{k=1}^N\psi_kI_{(t_{k-1},t_k]}(t),$$
where $0=t_0<t_1<\cds<t_{N-1}<t_N=T$ is a partition of $[0,T]$, and
each $A_k(\cd)$ is an $\dbF$-adapted process valued in
$\dbR^{n\times n}$,
$$A_k(s)\in\dbR^{n\times n}_{*+},\qq s\in[0,T],~\as$$
$$A_{k-1}(s)-A_k(s)\in\dbR^{n\times n}_+,\qq
s\in[0,T],~k=1,\cds,N,~\as,$$
each $\psi_k$ is an $\cF_T$-measurable random variable valued in
$\dbR^n$ such that
$$\psi_1\ge\psi_2\ge\cds\ge\psi_{N-1}\ge\psi_N\ge0,\qq\as$$
Let $(Y(\cd),Z(\cd\,,\cd))$ be the adapted solution to the BSVIE. On
$(t_{N-1},t_N]$, we have
$$Y(t)=\psi_N+\int_t^T\(A_N(s)Y(s)+B(s)Z(t,s)\)ds-\int_t^TZ(t,s)dW(s).$$
By uniqueness of BSDEs, we see that
$$(Y(s),Z(t,s))\equiv(Y_N(s),Z_N(s)),\qq\forall t_{N-1}<t\le s\le T,$$
with $(Y_N(\cd),Z_N(\cd))$ being the adapted solution to the
following BSDE:
$$Y_N(t)=\psi_N+\int_t^T\(A_N(s)Y_N(s)+B(s)Z_N(s)\)ds-\int_t^TZ_N(s)dW(s),\q t\in(t_{N-1},t_N].$$
Further, under our condition, by Proposition 2.4, we have
$$Y(t)\equiv Y_N(t)\ge0,\qq t\in(t_{N-1},t_N],~\as$$
In particular,
$$Y(t_{N-1}+0)=\psi_N+\int_{t_{N-1}}^T\(A_N(s)Y(s)+B(s)Z_N(s)\)ds
-\int_{t_{N-1}}^TZ_N(s)dW(s)\ge0,\qq\as$$
Next, for $t\in(t_{N-2},t_{N-1}]$, we have
$$\ba{ll}
\ns\ds
Y(t)=\psi_{N-1}+\int_t^T\(A_{N-1}(s)Y(s)+B(s)Z(t,s)\)ds-\int_t^TZ(t,s)dW(s)\\
\ns\ds\qq=\psi_{N-1}-\psi_N+Y(t_{N-1}+0)+\int_{t_{N-1}}^T\big[A_{N-1}(s)-A_N(s)\big]Y(s)ds\\
\ns\ds\qq\qq+\int_{t_{N-1}}^TB(s)\big[Z(t,s)-Z_N(s)\big]ds-\int_{t_{N-1}}^T\big[Z(t,s)
-Z_N(s)\big]dW(s)\\
\ns\ds\qq\qq+\int_t^{t_{N-1}}\(A_{N-1}(s)Y(s)+B(s)Z(t,s)\)ds-\int_t^{t_{N-1}}Z(t,s)dW(s).\ea$$
Let $(\wt Y_N(\cd),\wt Z_N(\cd))$ be the adapted solution to the
following BSDE:
$$\ba{ll}
\ns\ds\wt Y_N(\t)=\psi_{N-1}-\psi_N+Y(t_{N-1}+0)+\int_{t_{N-1}}^T\big[A_{N-1}(s)-A_N(s)\big]Y(s)ds\\
\ns\ds\qq\qq+\int_\t^TB(s)\wt Z_N(s)ds-\int_\t^T\wt
Z_N(s)dW(s),\qq\t\in(t_{N-1},T].\ea$$
Since
$$\psi_{N-1}-\psi_N+Y(t_{N-1}+0)+\int_{t_{N-1}}^T\big[A_{N-1}(s)-A_N(s)\big]Y(s)ds\ge0,$$
by our conditions, using Proposition 2.4, we have
$$\wt Y_N(\t)\ge0,\qq\t\in(t_{N-1},T],~\as$$
In particular,
$$\ba{ll}
\ns\ds\wt Y_N(t_{N-1}+0)=\psi_{N-1}-\psi_N+Y(t_{N-1}+0)+\int_{t_{N-1}}^T\big[A_{N-1}(s)-A_N(s)\big]Y(s)ds\\
\ns\ds\qq\qq\qq\qq+\int_{t_{N-1}}^TB(s)\wt
Z_N(s)ds-\int_{t_{N-1}}^T\wt Z_N(s)dW(s)\ge0,\qq\as\ea$$
On the other hand, by the uniqueness of adapted solutions to the
above BSVIEs, it is necessary that
$$Z(t,s)=
Z_N(s)+\wt Z_N(s),\qq(t,s)\in(t_{N-2},t_{N-1}]\times(t_{N-1},t_N].$$
Then $\wt Y_N(t_{N-1})$ is $\cF_{t_{N-1}}$-measurable, $s\mapsto
Z(t,s)$ is $\dbF$-adapted, and for $t\in(t_{N-2},t_{N-1}]$, and
$$Y(t)=\wt
Y_N(t_{N-1}+0)+\int_t^{t_{N-1}}\(A_{N-1}(s)Y(s)+B(s)Z(t,s)\)ds-\int_t^{t_{N-1}}Z(t,s)dW(s).$$
Next, we let $(Y_{N-1}(\cd),Z_{N-1}(\cd))$ be the adapted solution
to the following BSDE:
$$\ba{ll}
\ns\ds Y_{N-1}(t)=\wt
Y(t_{N-1}+0)+\int_t^{t_{N-1}}\(A_{N-1}(s)Y_{N-1}(s)+B(s)Z_{N-1}(s)\)ds\\
\ns\ds\qq\qq\qq-\int_t^{t_{N-1}}Z_{N-1}(s)dW(s),\qq
t\in[t_{N-2},t_{N-1}].\ea$$
By uniqueness of adapted solutions to BSDEs, we must have
$$(Y(s),Z(t,s))=(Y_{N-1}(s),Z_{N-1}(s)),\qq t\in(t_{N-2},t_{N-1}],~s\in[t,t_{N-1}].$$
Also, by $\wt Y(t_{N-1})\ge0$, we obtain
$$Y(t)\equiv Y_{N-1}(t)\ge0,\qq t\in(t_{N-2},t_{N-1}].$$
Therefore,
$$Y(t)\ge0,\qq t\in(t_{N-2},t_N],~\as$$
Then, by induction, we obtain
$$Y(t)\ge0,\qq t\in[0,T].$$
Finally, by approximation, we obtain the general case. \endpf

\ms

In the above proof, the condition that the coefficient $B(s)$ of
$Z(t,s)$ is independent of $t$ is crucial. It is desired if the
above remains true when $B(s)$ is replaced by $B(t,s)$.
Unfortunately, we do not have a confirmative answer at the
moment.

Having the above result, we now state a result for nonlinear case.

\ms

\bf Theorem 3.7. \sl Let
$g^i:\D\times\dbR^n\times\dbR^n\times\O\to\dbR^n$ satisfy
{\rm(BV1)}, and the following hold
\bel{gi}g^i(t,s,y,z)=h^i(t,s,y)+B(s)z,\qq(t,s,y,z)\in\D\times\dbR^n\times\dbR^n,\ee
for some $h^i:\D\times\dbR^n\times\O\to\dbR^n$ and $B(\cd)\in
L^\infty_\dbF(0,T;\dbR^{n\times n})$. Moreover,
\bel{3.23}\ba{ll}
\ns\ds h^1(t,s,y)-h^0(t,s,y)\ge
h^1(\t,s,y)-h^0(\t,s,y)\ge0,\\
\ns\ds\qq\qq\qq\qq\forall y\in\dbR^n,~0\le t\le\t\le s\le
T,~\as,\ea\ee
and for either $i=0$ or $i=1$, $y\mapsto h^i(t,s,y)$ is
differentiable with
\bel{3.19_1}h^i_y(t,s,y)\in\dbR^{n\times n}_{*+},\qq
h^i_y(t,s,y)-h^i_y(\t,s,y)\in\dbR^{n\times n}_+,\qq0\le t\le\t\le
s\le T,~y\in\dbR^n,~\as\ee
Then for any $\psi^i(\cd)\in C_{\cF_T}([0,T];L^2(\O;\dbR^n))$ with
\bel{}\psi^1(t)-\psi^0(t)\ge\psi^1(\t)-\psi^0(\t)\ge0,\qq0\le
t\le\t\le T,~\as,\ee
the corresponding adapted solutions $(Y^i(\cd),Z^i(\cd\,,\cd))$ of
BSVIEs $(\ref{BSVIE5})$ satisfy
$$Y^1(t)\ge Y^0(t),\qq t\in[0,T],~\as$$

\ms

\rm

\it Proof. \rm Suppose that $y\mapsto h^0(t,s,y)$ is differentiable
and (\ref{3.19_1}) holds for $i=0$. Then we have
$$\ba{ll}
\ns\ds
Y^1(t)-Y^0(t)=\psi^1(t)-\psi^0(t)+\int_t^T\[h^1(t,s,Y^1(s))-h^0(t,s,Y^1(s))\]ds\\
\ns\ds\qq\qq+\int_t^T\[A(t,s)\(Y^1(s)-Y^0(s)\)+B(s)\(Z^1(t,s)-Z^0(t,s)\)\]ds\\
\ns\ds\qq\qq-\int_t^T\(Z^1(t,s) -Z^0(t,s)\)dW(s),\ea$$
where
$$A(t,s)=\int_0^1
h^0_y(t,s,Y^0(s)+\b[Y^1(s)-Y^0(s)])d\b,\qq(t,s)\in\D.$$
Then our conclusion follows from Theorem 3.6. \endpf

\ms

Note that in the above theorem, we have not assumed any sort of
nondecreasing conditions on $y\mapsto g^i(t,s,y,z)$. Also, when
$h^i(t,s,y)$ are independent of $t$, condition (\ref{3.18}) is
reduced to
$$h^1(s,y)\ge h^0(s,y),\qq(s,y)\in[0,T]\times\dbR^n,~\as,$$
and condition (\ref{3.19_1}) is automatically true. Finally, if
$y\mapsto h^i(t,s,y)$ is just Lipschitz and not necessarily
differentiable, we may modify condition (\ref{3.19_1}) in a proper
way so that the same conclusion remains. On the other hand, we have
seen that our result does not fully recover the comparison theorem
for general nonlinear $n$-dimensional BSDEs. At the moment, this is
the best that we can do.

\ms

\subsection{Comparison theorem for adapted M-solutions.}

In this subsection, we discuss the following type BSVIEs:
\bel{BSVIE6}\ba{ll}
\ns\ds
Y(t)=\psi(t)+\int_t^Tg(t,s,Y(s),Z(s,t))ds-\int_t^TZ(t,s)dW(s),\qq
t\in[0,T].\ea\ee
Note that since the generator $g(\cd)$ depends on $Z(s,t)$, the
notion of adapted solution in $\mathcal{H}^p_{\Delta}[0,T]$ will not
be enough. Therefore, we adopt the notion of adapted M-solution to
the above BSVIE (\cite{Yong 2008}). More precisely, an adapted
M-solution is an adapted solution $(Y(\cd),Z(\cd\,,\cd))$ which
belongs to $\cM^p[0,T]$. The following is a standard assumption for
the BSVIE (\ref{BSVIE6}).

\ms

{\bf(BV2)} For $i=0,1$, the maps
$g^i:\D\times\dbR^n\times\dbR^n\times\O\to\dbR^n$ is measurable,
$s\mapsto g^i(t,s,y,\z)$ is $\dbF$-progressively measurable,
$(y,\z)\mapsto g^i(t,s,y,\z)$ is uniformly Lipschitz, $(t,s)\mapsto
g^i(t,s,0,0)$ is uniformly bounded.

\ms

By \cite{Yong 2008}, we know that under (BV2), for any $\psi(\cd)\in
C_\dbF([0,T];L^2(\O;\dbR^n))$, (\ref{BSVIE6}) admits a unique
adapted M-solution $(Y(\cd),Z(\cd\,,\cd))$. We will use a dual
principle (\cite{Yong 2008}) to prove the comparison theorem for
adapted M-solution. The results of this subsection also corrects
relevant ones in \cite{Yong 2006, Yong 2007}. Before going further,
let us look at a simple example.

\ms

\bf Example 3.8. \rm Consider the following one-dimensional BSVIE:
\bel{3.2.2}\ba{ll}
\ns\ds
Y(t)=\psi(t)+\int_t^T\frac{2T-t}{2T-s}Z(s,t)ds-\int_t^TZ(t,s)dW(s),\qq
t\in[0,T].\ea\ee
We introduce the following FSVIE:
\bel{3.2.3}X(t)=1+\int_0^t\frac{2T-s}{2T-t}X(s)dW(s),\qq
t\in[0,T],\ee
which is the equation in Example 2.7, and
$$\mathbb{P}\{\omega;X(t,\omega)<0\}>0,\quad \forall t\in[0,T].$$
Hence by taking
$$\psi(t)=I_{A}(t,\omega),\quad A=\{(t,\omega);X(t,\omega)<0\},\ t\in[0,T],$$
by the duality principle (\cite{Yong 2008}), we have
$$\dbE\int_0^TY(t)dt=\dbE\int_0^TX(t)I_{A}(t,\omega)dt<0,$$
which means that $Y(\cd)\ge0$ on $[0,T]$ could not be true, although
$\psi(\cd)\ge0$.

\ms

The above example shows that comparison theorem may fail for linear
BSVIEs if in the generator, the coefficient of $Z(s,t)$ depends both
on $t$ and $s$. The above example suggests us that if linear BSVIEs
are considered for comparison of adapted M-solutions, the following
should be a proper form:
\bel{3.24}\ba{ll}
\ns\ds
Y(t)=\psi(t)+\int_t^T\(A(t,s)Y(s)+C(t)Z(s,t)\)ds-\int_t^TZ(t,s)dW(s),\qq
t\in[0,T].\ea\ee
Note that $Z(t,s)$ does not appear in the drift term, and the
coefficient $C(t)$ of $Z(s,t)$ is independent of $s$. For such an
equation, we have the following result, which is comparable with
Theorem 3.6.

\ms

\bf Theorem 3.9. \sl Let $A:\D\times\O\to\dbR^{n\times n}$ and
$C:[0,T]\times\O\to\dbR^{n\times n}$ be uniformly bounded, with
$C(\cd)$ being $\dbF$-progressively measurable, for each
$t\in[0,T]$, $s\mapsto A(t,s)$ being $\dbF$-progressively
measurable, and for each $s\in[0,T]$, $t\mapsto A(s,t)$ is
continuous. Further,
\bel{}A(t,s)\in\dbR_{*+}^{n\times n},\qq(t,s)\in\D,~\as,\ee
\bel{}A(s,\t)-A(s,t)\in\dbR^{n\times n}_+,\qq\forall\,s\le t\le\t\le
T,~s\in[0,T],~\as,\ee
\bel{}C(t)\in\dbR_d^{n\times n},\qq\ae t\in[0,T],~\as\ee
Then the adapted M-solution $(Y(\cd),Z(\cd\,,\cd))$ of linear BSVIE
$(\ref{3.24})$ with $\psi(\cd)\in C_{\cF_T}(0,T;L^2(\O;\dbR^n))$,
$\psi(\cd)\ge0$ satisfies
\bel{3.33}\dbE_t\int_t^TY(s)ds\ge0,\qq\forall t\in[0,T],~\as\ee

\ms

\it Proof. \rm Pick any $\eta(\cd)\in L^2_\dbF(0,T;\dbR^n)$ with
$\eta(\cd)\ge0$, consider the following linear FSVIE:
\bel{}\ba{ll}
\ns\ds X(t)=\f(t)+\int_0^tA(s,t)^TX(s)ds+\int_0^t
C(s)^TX(s)dW(s),\qq t\in[0,T],\ea\ee
with
$$\f(t)=\int_0^t\eta(s)ds,\qq t\in[0,T].$$
By our conditions on $A(\cd\,,\cd)$ and $C(\cd)$, using Proposition
2.7, we have
$$X(t)\ge0,\qq t\in[0,T],~\as$$
Then by duality theorem (\cite{Yong 2008}), one obtains
$$\ba{ll}
\ns\ds0\le\dbE\int_0^T\lan\psi(t),X(t)\ran
dt=\dbE\int_0^T\lan\f(t),Y(t)\ran dt\\
\ns\ds\q=\dbE\int_0^T\int_0^t\lan\eta(s),Y(t)\ran
dsdt=\dbE\int_0^T\lan\eta(s),\int_s^TY(t)dt\ran ds.\ea$$
Thus (\ref{3.33}) follows since $\eta(\cd)$ is arbitrary. \endpf

\ms

Different from Theorem 3.6, in the above, we do not need the
monotonicity of $t\mapsto\psi(t)$, and the conclusion (\ref{3.33})
is weaker than (\ref{3.22}).

\ms

Having the above result, we are able to get a comparison theorem for
the following nonlinear BSVIEs ($i=0,1$)
\bel{5.34}Y^i(t)=\psi^i(t)+\int_t^T\(h^i(t,s,Y^i(s))+
C(t)Z^i(s,t)\)ds-\int_t^TZ^i(t,s)dW(s),\q t\in[0,T].\ee
More precisely, the following theorem holds.

\ms

\bf Theorem 3.10. \sl Let
$g^i:\D\times\dbR^n\times\dbR^n\times\O\to\dbR^n$ satisfy {\rm(BV2)}
and the following hold
\bel{gi}g^i(t,s,y,\z)=h^i(t,s,y)+C(t)\z,\qq(t,s,y,\z)\in\D\times\dbR^n\times\dbR^n,\ee
for some $h^i:\D\times\dbR^n\times\O\to\dbR^n$ and $C(\cd)\in
L^\infty_\dbF(0,T;\dbR^{n\times n})$. Moreover,
\bel{3.36}\ba{ll}
\ns\ds h^1(t,s,y)-h^0(t,s,y)\ge h^1(\t,s,y)-h^0(\t,s,y)\ge0,\qq
\forall y\in\dbR^n,~0\le t\le\t\le s\le T,~\as,\ea\ee
and for either $i=0$ or $i=1$, $y\mapsto h^i(t,s,y)$ is
differentiable with
\bel{3.19}h^i_y(t,s,y)\in\dbR^{n\times n}_{*+},\qq
h^i_y(t,s,y)-h^i_y(\t,s,y)\in\dbR^{n\times n}_+,\qq0\le t\le\t\le
s\le T,~y\in\dbR^n,~\as\ee
Then for any $\psi^i(\cd)\in C_{\cF_T}([0,T];L^2(\O;\dbR^n))$ with
\bel{}\psi^1(t)-\psi^0(t)\ge\psi^1(\t)-\psi^0(\t)\ge0,\qq0\le
t\le\t\le T,~\as,\ee
the corresponding adapted solutions $(Y^i(\cd),Z^i(\cd\,,\cd))$ of
BSVIEs $(\ref{BSVIE6})$ satisfy
\bel{3.40}\dbE_t\int_t^TY^1(s)ds\ge\dbE_t\int_t^TY^0(s)ds,\qq
t\in[0,T],~\as\ee

\it Proof. \rm Observe the following:
$$\ba{ll}
\ns\ds
Y^1(t)-Y^0(t)=\psi^1(t)-\psi^0(t)+\int_t^T\[h^1(t,s,Y^1(s))-h^0(t,s,Y^1(s))\]ds\\
\ns\ds\qq\qq\qq\qq+\int_t^T\[A(t,s)\(Y^1(s)-Y^0(s)\)
+C(t)\(Z^1(s,t)-Z^0(s,t)\)\]ds\\
\ms\ds\qq\qq\qq\qq-\int_t^T\(Z^1(t,s)-Z^0(t,s)\)dW(s),\ea$$
where
$$A(t,s)=\int_0^1h^0_y\big(t,s,Y^0(s)+\b[Y^1(s)-y^0(s)]\big)d\b.$$
Then under our conditions, we have the comparison (\ref{3.40}).
\endpf

\subsection{Other type solutions to BSVIEs.}

We now look at the following general BSVIE:
\bel{BSVIE7}Y(t)=\psi(t)+\int_t^Tg(t,s,Y(s),Z(t,s),Z(s,t))ds-\int_t^TZ(t,s)dW(s),\q
t\in[0,T].\ee
According to \cite{Yong 2008}, due to the appearance of $Z(s,t)$,
there are infinite adapted solutions for (\ref{BSVIE7}), and under
proper conditions, (\ref{BSVIE7}) admits a unique adapted M-solution
$(Y(\cd),Z(\cd\,,\cd))$. On the other hand, it is possible to define
other types of solutions.

\ms

Recall that in the mean-variance problem, the precommitted solution
is widely used when the objective function is
$$\dbE X(T)-{\g\over2}\dbE\(X(T)-\dbE X(T)\)^2,$$
with $X(T)$ being the terminal wealth. Recently people started to
study the dynamic version of
$$\dbE_tX(T)-{\g\over2}\dbE_t\[X(T)-\dbE_tX(T)\]^2,\qq t\in[0,T],$$
and proposed the time consistent solution, see
for example \cite{BMZ2012}. On the other hand, mean-field BSDE of
\bel{MF-BSVIE}Y(t)=\xi+\int_t^Tg(s,Y(s),Z(s),\dbE Y(s),\dbE
Z(s))ds-\int_t^TZ(s)dW(s),\qq t\in[0,T],\ee
was introduced and studied in \cite{BLP2009}. A dynamic version of
(\ref{MF-BSVIE}) should be the following:
\bel{3.44}Y(t)=\psi(t)+\int_t^Tg(s,Y(s),Z(t,s),\dbE_tY(s),
\dbE_tZ(t,s))ds-\int_t^TZ(t,s)dW(s),\q t\in[0,T].\ee
Inspiring by the above, we introduce the following definition.

\ms

\bf Definition 3.11. \rm Let
$h:\D\times\dbR^n\times\dbR^n\times\O\to\dbR^n$ satisfy (BV1).
Moreover, $t\mapsto h(t,s,y,z)$ is $\mathcal{F}_t$-measurable for
given $(s,y,z)\in[0,T]\times\dbR^n\times\dbR^n$. A pair
$(Y(\cd),Z(\cd\,,\cd))\in\cH^2[0,T]$ is called a {\it conditional
$h$-solution} for BSVIE (\ref{BSVIE7}) if (\ref{BSVIE7}) is
satisfied in the It\^o sense and
\bel{3.45}Z(s,t)=h(t,s,\dbE_tY(s),\dbE_tZ(t,s)),\qq
(t,s)\in\D,~\as\ee

\ms

It is clear that for BSVIE (\ref{BSVIE7}), if we are talking about
conditional $h$-solution, it amounts to studying (the usual) adapted
solution for BSVIE (\ref{3.44}) of mean-field type. By using a
similar method in \cite{Yong 2008} or \cite{Wang-Shi-Yong 2011}, we
can establish the existence and uniqueness of adapted solution to
(\ref{3.44}). Then following the ideas contained in the previous
subsections, we are able to discuss comparison of adapted solutions
for such kind of equations. We prefer not to get into details here.

\ms

\section{Concluding Remarks.}

For BSIVEs of form
$$Y(t)=\psi(t)+\int_t^Tg(t,s,Y(s),Z(t,s))ds-\int_t^TZ(t,s)dW(s),\qq t\in[0,T],$$
we have established a general comparison theorem (Theorem 3.2) for
the adapted solutions when the tube
$$\Big\{[g^0(t,s,y,z),g^1(t,s,y,z)]\big|(t,s,y,z)\in\D\times\dbR^n\times\dbR^n\Big\}$$
admits a selection $\bar g(t,s,y,z)$ which is nondecreasing in $y$,
plus some additional conditions. Examples 3.3, 3.4, 3.5, and 3.7
tell us that when the above condition is not assumed, the situation
becomes very complicated. At the moment, if the above monotonicity
condition is not assumed, we can only prove a comparison theorem for
the following restricted form of BSVIEs (see Theorem 3.6):
$$Y(t)=\psi(t)+\int_t^T\(h(t,s,Y(s))+B(s)Z(t,s)\)ds-\int_t^TZ(t,s)dW(s),\qq
t\in[0,T].$$
Further, if the generator depends on $Z(s,t)$, we need to use
duality principle to prove a proper comparison theorem. Due to this,
at the moment, the BSVIEs that we can treat is the following type:
$$Y(t)=\psi(t)+\int_t^T\(h(t,s,Y(s))+C(t)Z(s,t)\)ds-\int_t^TZ(t,s)dW(s),\qq
t\in[0,T].$$
Moreover, if $(Y^i(\cd),Z^i(\cd\,,\cd))$ ($i=0,1$) are adapted
M-solutions to the BSVIEs of the above form, instead of
$$Y^0(t)\le Y^1(t),\qq t\in[0,T],~\as,$$
(under suitable conditions, see Theorem 3.11), we only have a weaker
form of comparison:
$$\dbE_t\[\int_t^TY^0(s)ds\]\le\dbE_t\[\int_t^TY^1(s)ds\],\qq t\in[0,T],~\as$$
Theorems 3.2, 3.6, and 3.10 correct the relevant result presented in
\cite{Yong 2006,Yong 2007}. Finally, the problem of comparison for
the adapted M-solutions to the following general type BSVIEs:
$$Y(t)=\psi(t)+\int_t^Tg(t,s,Y(s),Z(t,s),Z(s,t))ds-\int_t^TZ(t,s)dW(s),\qq
t\in[0,T],$$
is widely open at the moment. We hope that some further results
could be addressed in our future publications.

\end{document}